\documentclass[preprint,12pt]{elsarticle}

\usepackage{ulem}
\usepackage{indentfirst}
\usepackage{bm}
\usepackage{graphicx,amsmath,amsfonts,color,mathrsfs, amsmath,amsthm}
\usepackage{amssymb}
\biboptions{numbers,sort&compress}
\usepackage{epsfig}
\begin{document}
\newcommand{\BOX}{\hfill $\Box$}
\newcommand{\eqn}[1]{(\ref{eqn:#1})}
\newcommand{\lem}[1]{Lemma~\ref{lem:#1}}
\newcommand{\cor}[1]{Corollary~\ref{cor:#1}}
\newcommand{\thr}[1]{Theorem~\ref{thr:#1}}
\newcommand{\con}[1]{Conjecture~\ref{con:#1}}
\newcommand{\pro}[1]{Proposition~\ref{pro:#1}}
\newcommand{\que}[1]{Question~\ref{que:#1}}
\newcommand{\exa}[1]{Example~\ref{exa:#1}}
\newcommand{\den}[1]{Definition~\ref{den:#1}}
\newcommand{\ass}[1]{Assumption~\ref{ass:#1}}
\newcommand{\rem}[1]{Remark~\ref{rem:#1}}
\newcommand{\fig}[1]{Figure~\ref{fig:#1}}

\numberwithin{equation}{section}

\newtheorem{theorem}{Theorem}[section]
\newtheorem{lemma}{Lemma}[section]
\newtheorem{assumption}{Assumption}[section]
\newtheorem{proposition}{Proposition}[section]
\newtheorem{conjecture}{Conjecture}[section]
\newtheorem{remark}{Remark}[section]
\newtheorem{question}{Question}[section]
\newtheorem{problem}{Problem}[section]
\newtheorem{example}{Example}[section]
\newtheorem{result}{Result}[section]
\newtheorem{corollary}{Corollary}[section]
\newtheorem{definition}{Definition}[section]

\begin{frontmatter}

%% Title, authors and addresses

%% use the tnoteref command within \title for footnotes;
%% use the tnotetext command for theassociated footnote;
%% use the fnref command within \author or \address for footnotes;
%% use the fntext command for theassociated footnote;
%% use the corref command within \author for corresponding author footnotes;
%% use the cortext command for theassociated footnote;
%% use the ead command for the email address,
%% and the form \ead[url] for the home page:
%% \title{Title\tnoteref{label1}}
%% \tnotetext[label1]{}
%% \author{Name\corref{cor1}\fnref{label2}}
%% \ead{email address}
%% \ead[url]{home page}
%% \fntext[label2]{}
%% \cortext[cor1]{}
%% \address{Address\fnref{label3}}
%% \fntext[label3]{}

\title{Analysis of the equilibrium strategies in the Geo/Geo/1 queue with multiple working vacations}

%% use optional labels to link authors explicitly to addresses:
%% \author[label1,label2]{}
%% \address[label1]{}
%% \address[label2]{}

\author{Bixuan Yang}
\ead{bixuanyang@126.com}
\author{Zhenting Hou}
\ead{zthou@csu.edu.cn}
\author{Jinbiao Wu \corref{cor1}}
\ead{wujinbiao@csu.edu.cn}
\author{Zaiming Liu}
\ead{math\_lzm@csu.edu.cn}
 \cortext[cor1]{Corresponding author}
\address{School of Mathematics and Statistics, Central South University, Changsha 410083, China}

\begin{abstract}
This paper studies the equilibrium behavior of customers in the Geo/Geo/1 queueing system with multiple working vacations. The arriving customers decide whether to join or to balk the queueing systems based on the information of the queue length and the states of the server. In observable queues, partially observable queues and unobservable queues three cases we obtain the equilibrium balking strategies based on the reward-cost structure and socially optimal strategies for all customers. Furthermore, we present some numerical experiments to illustrate the effect of the information level on the equilibrium behavior and to compare the customers' equilibrium and socially optimal strategies.

\end{abstract}

\begin{keyword}
%% keywords here, in the form: keyword \sep keyword

%% PACS codes here, in the form: \PACS code \sep code

%% MSC codes here, in the form: \MSC code \sep code
%% or \MSC[2008] code \sep code (2000 is the default)
Geo/Geo/1 queue\sep Multiple working vacations\sep Equilibrium balking strategies\sep Socially optimal strategies
\end{keyword}

\end{frontmatter}

%% \linenumbers

%% main text

\section{Introduction}
\label{1}

%% The Appendices part is started with the command \appendix;
%% appendix sections are then done as normal sections
%% \appendix

%% \section{}
%% \label{}

%% If you have bibdatabase file and want bibtex to generate the
%% bibitems, please use
%%
%%  \bibliographystyle{elsarticle-harv}
%%  \bibliography{<your bibdatabase>}

%% else use the following coding to input the bibitems directly in the
%% TeX file.
Due to important applications in the fields of computer networks, communications systems and production management, economic queueing systems have received comprehensive attention. Recently, based on different queueing models to study the customers' decisions that whether to join or to balk the queueing systems have became a new hot spot. Such an economic analysis of queueing systems was pioneered by Naor \cite{1}, who studied the equilibrium and socially optimal strategies in an $M/M/1$ queue with a simple linear reward-cost structure. On Naor's model, arriving customers were informed about the queue length before they made decisions. Edelson and Hildeband \cite{2} considered the unobservable case in which the customers made their decisions without being informed about the state of the queue. Subsequently, Naor's model and results had been extended in several literatures, see e.g. \cite{3, 4, 5}. Larsen \cite{6} promoted a more generalized model assuming that customers' service values were a set of random variables instead of fixed constants. Chen and Frank \cite{7} spreaded Naor's model assuming that both the customers and the server maximize their expected discounted utility using a common discount rate. Mandelbaum and Shimkin \cite{8}, Shimkin and Mandelbaum \cite{9} respectively discussed the equilibrium conditions when the expenditure functions were linear and non-linear in non-visual systems. The fundamental results on this subject in both the observable and unobservable queueing systems can be found in the comprehensive monographs of Hassin and Haviv \cite{10}.

Discrete-time queueing systems with vacations have been widely researched by a lot of investigators because of their extensively applications in manufacturing and telecommunication systems. An excellent and complete study on discrete-time queueing systems with vacations had been presented by Takagi \cite{11}. Zhang and Tian \cite{12} presented the detailed analysis on the $Geo/G/1$ queue with multiple adaptive vacations. Moreover, Goswami and Mund \cite{13} analyzed a finite-buffer renewal input single server discrete-time $GI/Geo/1/N$ queue with
multiple working vacations. Recently, Vijaya Laxmi and Seleshi \cite{14} investigated a discrete-time renewal input queue with change
over times and Bernoulli schedule vacation interruption under batch service $(a, c, b)$ policy.

The study about the equilibrium customer behavior in vacation queue models, the first was presented by Burnetas and Economou \cite{15}, who explored both the observable and unobservable conditions in a single server Markovian queue with setup times. Then, Economou and Kanta \cite{16} examined the equilibrium balking strategies and pricing for the single server Markovian queue with compartmented waiting space. Guo et al. \cite{17} considered customer equilibrium and socially optimal strategies to join a queue with only partial information on the service time distribution such as moments and the range. Liu et al. \cite{18} researched the equilibrium threshold strategies in observable queueing systems under single vacation policy. Sun et al. \cite{19} further studied the customers' equilibrium and socially optimal joining-balking behavior in the $M/M/1$ queue with multiple working vacations. However, to the best of the authors' knowledge, researches for the equilibrium strategies in the discrete-time queueing systems with multiple working vacations haven't been given.

In this paper, we investigate the equilibrium joining/balking behavior of the customers in the discrete-time $Geo/Geo/1$ queueing system with multiple working vacations. The server takes the original work at the lower rate rather than completely stopping during the working vacation period. We will explore the equilibrium strategies for customers' decisions that whether to join or to balk the queueing systems based on different information upon arrival. Three cases will be considered: (1) Observable case: arriving customers are informed about both the states of the server and the number of customers in the system; (2) Partially observable case: arriving customers are informed only about the states of the server; (3) Unobservable case: arriving customers are informed about neither the states of the server nor the number of customers in the system. In the three cases we study the equilibrium strategies based on the reward-cost structure and socially optimal strategies for all customers. Subsequently, we present some numerical experiments to research the effect of the information level on the equilibrium behavior and to compare the customers' equilibrium and socially optimal strategies.

This paper is organized as follows. In Section 2, we describe the queueing model. In Sections 3, 4, and 5, we respectively study the observable queues, the partially observable queues and the unobservable queues by using the equilibrium threshold strategies and present some numerical examples to explain the effect of the information level on customers' behavior and to compare the customers' equilibrium and socially optimal strategies. In Section 6, the conclusions are given.

\section{Description of the model}
\label{2}

Geo/Geo/1 queueing system with multiple working vacations can be described as follows. Throughout this paper, for any real number $x\in[0,1]$, we denote $\bar{x}=1-x$.\\
(1) Potential customer arrivals occur at the end of slot $(n^{-},n)$. The inter-arrival times are independent and identically distributed sequences which follow a  geometric distribution with rate $p$, $0<p<1$,
\[
P(T=k)=p\bar{p}^{k-1},\ \ k\geq1.
\]
(2) The beginning and ending of the service take place at slot division point $t=n$, $n=0,1,2,...$. The regular service times are independent each other and geometrically distributed with rate $\mu_{b}$, $0<\mu_{b}<1$,
\[
P(S_{b}=k)=\mu_{b}\bar{\mu}_{b}^{k-1},\ \ k\geq1.
\]
In a working vacation period, the service times are independent each other and geometrically distributed with rate $\mu_{\nu}$, $0<\mu_{\nu}<1$,
\[
P(S_{\nu}=k)=\mu_{\nu}\bar{\mu}_{\nu}^{k-1},\ \ k\geq1.
\]
(3) When the queue becomes empty, the server enters a working vacation period. The working vacation time $V$ follows a geometric distribution with parameter $\theta$, $0<\theta<1$,
\[
P(V=k)=\theta\bar{\theta}^{k-1},\ \ k\geq1.
\]
During a working vacation arriving customers are served according to arrival order by the rate $\mu_{\nu}$. When a working vacation ends, if there are customers in the queue, the server switches the service rate from $\mu_{\nu}$ to $\mu_{b}$, and a regular busy period begins. Otherwise, the server keeps on another working vacation. Suppose that the beginning and ending of vacations occur at the end of slot $(n^{-},n)$.\\
(4) Assume that inter-arrival times, service times and working vacation times are mutually independent. The queueing system follows the First-Come-First-Served (FCFS) service discipline.

Let $L_{n}^{+}$ be the number of customers in the system at time $n^{+}$. According to the above description, a customer who completes service and leaves at $(n,n^{+})$ no longer be included in $L_{n}^{+}$. Define
$$
     J_{n}=\left\{\begin{array}{ll}
                    0, \mbox{the system is in a working vacation period at time}~ n^{+},\\
                    1, \mbox{the system is in a regular busy period at time}~ n^{+}.
                  \end{array}
\right.
$$
It is easy to get that $\{(L_{n}^{+},J_{n}),n\geqslant0\}$ is a Markov chain with state space
\begin{center}
$\Omega=\{(0,0)\}\bigcup\{(k,j):k\geqslant1,~j=0,1\}$.
\end{center}

We investigate customer equilibrium strategies for joining/balking. We distinguish three cases with respect to the level of information available to customers at their arrival instants, before their decisions are made: the observable case means customers can observe both $J_{n}$ and $L_{n}^{+}$; the partially observable case means customers can observe $J_{n}$ but not $L_{n}^{+}$; the unobservable case means customers can observe neither $J_{n}$ nor $L_{n}^{+}$.

Our interest is the behavior of the customers when they decide whether to join or to balk upon their arrival. To model the decision process, we assume that every customer receives a reward of $R$ units after service completion. It may reflect a customer's satisfaction and the added value of being served. Moreover, there exists a waiting cost of $C$ units per time unit that a customer remains in the system (in queue and in service area). In this paper, we use a linear cost function, then a customer's expected net benefit after service completion, denoted by $U$, is $U=R-CE[W]$, where $E[W]$ represents the customer's mean sojourn time in the queueing system. We can get that if the customer selects balking, $U=0$.

\section{The observable queues}
\label{3}

We first consider the observable case in which the arriving customers can observe both the state $J_{n}$ of the server and the number of customers $L_{n}^{+}$ in the system. We need to use equilibrium strategies of thresholds type. In the observable queues, a pure threshold strategy is specified by a pair $(n_{e}(0),n_{e}(1))$ and has the form `observe $(L_{n}^{+},J_{n})$ at arrival instant; enter if $L_{n}^{+}\leq n_{e}(J_{n})$ and balk otherwise'.

As for the mean sojourn time of a customer who joins the system at state $(n,1),n\geq1$, we need to consider as follows. Assume that a new arrival occurs at $(n^{-},n)$, it is possible that a service is ending. So the sojourn time equals $n$ service times with probability $\mu_{b}$ and equals $n+1$ service times with probability $\bar{\mu}_{b}$. We get the probability generating function (PGF) of the sojourn time of a customer who joins the system at state $(n,1)$, denoted by $\widetilde{W}_{1}^{*}(z)$.
\begin{equation*}
\begin{aligned}
\widetilde{W}_{1}^{*}(z)=\mu_{b}\left(\frac{\mu_{b}z}{1-\bar{\mu}_{b}z}\right)^{n}+\bar{\mu}_{b}\left(\frac{\mu_{b}z}{1-\bar{\mu}_{b}z}\right)^{n+1}
=\frac{\mu_{b}}{1-\bar{\mu}_{b}z}\left(\frac{\mu_{b}z}{1-\bar{\mu}_{b}z}\right)^{n}.
\end{aligned}
\end{equation*}
Then the mean sojourn time of a customer who joins the system at state $J=1$ is
\begin{equation*}
E[W(1)]=\widetilde{W}_{1}^{*'}(z)|_{z=1}=\frac{n+1}{\mu_{b}}-1.
\end{equation*}

As for the mean sojourn time of a customer who joins the system at state $(n,0),n\geq0$, we need to consider two situations. Denote $S_{\nu}^{(j)}$ the sum of $j$ service times $S_{\nu}$, $S_{\nu}^{(n)}(z|S_{\nu}^{(n)}\leq V)$ the PGF of $S_{\nu}^{(n)}$ under the condition $S_{\nu}^{(n)}\leq V$ and $V(z|S_{\nu}^{(j)}\leq V <S_{\nu}^{(j+1)})$ the PGF of $V$ under the condition $S_{\nu}^{(j)}\leq V <S_{\nu}^{(j+1)}$. One situation is that when a new arrival occurs at $(n^{-},n)$, there is a customer leaving the system after service completion at the instant $t=n$. The probability of this event is equal to $\mu_{\nu}$. Then there are two cases. Case 1: if there are $j$ customer service completions when the working vacation ends, $ 0 \leqslant j\leqslant n-1$, the sojourn time is the sum of the residual working vacation time under the condition $S_{\nu}^{(j)}\leq V <S_{\nu}^{(j+1)}$ plus $n-j$ service times with rate $\mu_{b}$. Case 2: if at least $n$ customers are served when the working vacation ends, the sojourn time is $n$ service times with rate $\mu_{\nu}$ under the condition $V\geqslant S_{\nu}^{(n)}$. Another situation is that an arrival occurs at $(n^{-},n)$, no customer is about to leave at the instant $t=n$. The probability of this event is equal to $\bar{\mu}_{\nu}$. There are still two cases and are similar to the above. So we get the PGF of the sojourn time of a customer who joins the system at state $(n,0)$, denoted by $\widetilde{W}_{0}^{*}(z)$.
\begin{equation*}
\begin{aligned}
\widetilde{W}_{0}^{*}(z)&=\mu_{\nu}\left[P(S_{\nu}^{(n)}\leq V)S_{\nu}^{(n)}(z|S_{\nu}^{(n)}\leq V)\right.\\
&~~~\left.+\sum_{j=0}^{n-1}P(S_{\nu}^{(j)}\leq V <S_{\nu}^{(j+1)})V(z|S_{\nu}^{(j)}\leq V<S_{\nu}^{(j+1)})\left(\frac{\mu_{b}}{1-\bar{\mu}_{b}z}\right)^{n-j}\right]\\
&~~~+\bar{\mu}_{\nu}\left[P(S_{\nu}^{(n+1)}\leq V)S_{\nu}^{(n+1)}(z|S_{\nu}^{(n+1)}\leq V)\right.\\
&~~~\left.+\sum_{j=0}^{n}P(S_{\nu}^{(j)}\leq V <S_{\nu}^{(j+1)})V(z|S_{\nu}^{(j)}\leq V <S_{\nu}^{(j+1)})\left(\frac{\mu_{b}}{1-\bar{\mu}_{b}z}\right)^{n+1-j}\right]\\
&~~~+P(S_{\nu}^{(1)}\leq V)S_{\nu}^{(1)}(z|S_{\nu}^{(1)}\leq V)+P(S_{\nu}^{(0)}\leq V <S_{\nu}^{(1)})V(z|S_{\nu}^{(0)}\leq V <S_{\nu}^{(1)})\frac{\mu_{b}}{1-\bar{\mu}_{b}z}\\
&=\mu_{\nu}\left[\left(\frac{\mu_{\nu}\bar{\theta}z}{1-\bar{\mu}_{\nu}\bar{\theta}z}\right)^{n}
+\sum_{j=0}^{n-1}\frac{\theta}{1-\bar{\mu}_{\nu}\bar{\theta}z}\left(\frac{\mu_{\nu}\bar{\theta}z}{1-\bar{\mu}_{\nu}\bar{\theta}z}\right)^{j}\left(\frac{\mu_{b}z}{1-\bar{\mu}_{b}z}\right)^{n-j}\right]\\
&~~~+\bar{\mu}_{\nu}\left[\left(\frac{\mu_{\nu}\bar{\theta}z}{1-\bar{\mu}_{\nu}\bar{\theta}z}\right)^{n+1}
+\sum_{j=0}^{n}\frac{\theta}{1-\bar{\mu}_{\nu}\bar{\theta}z}\left(\frac{\mu_{\nu}\bar{\theta}z}{1-\bar{\mu}_{\nu}\bar{\theta}z}\right)^{j}\left(\frac{\mu_{b}z}{1-\bar{\mu}_{b}z}\right)^{n+1-j}\right]\\
&~~~+\frac{\mu_{\nu}\bar{\theta}z}{1-\bar{\mu}_{\nu}\bar{\theta}z}+\frac{\theta}{1-\bar{\mu}_{\nu}\bar{\theta}z}\frac{\mu_{b}z}{1-\bar{\mu}_{b}z}\\
&=\frac{\mu_{\nu}}{1-\bar{\mu}_{\nu}\bar{\theta}z}\left(\frac{\mu_{\nu}\bar{\theta}z}{1-\bar{\mu}_{\nu}\bar{\theta}z}\right)^{n}\left(1-\frac{\theta\mu_{b}}{\mu_{b}(1-\bar{\mu}_{\nu}\bar{\theta}z)-\mu_{\nu}\bar{\theta}(1-\bar{\mu}_{b}z)}\right)\\
&~~~+\frac{\theta\mu_{b}}{\mu_{b}(1-\bar{\mu}_{\nu}\bar{\theta}z)-\mu_{\nu}\bar{\theta}(1-\bar{\mu}_{b}z)}\left(\frac{\mu_{b}z}{1-\bar{\mu}_{b}z}\right)^{n}\left(\mu_{\nu}+\frac{\bar{\mu}_{\nu}\mu_{b}z}{1-\bar{\mu}_{b}z}\right)\\
&~~~+\frac{\mu_{\nu}\bar{\theta}z}{1-\bar{\mu}_{\nu}\bar{\theta}z}+\frac{\theta}{1-\bar{\mu}_{\nu}\bar{\theta}z}\frac{\mu_{b}z}{1-\bar{\mu}_{b}z}.
\end{aligned}
\end{equation*}
Then the mean sojourn time of a customer who joins the system at state $J=0$ is
\begin{equation*}
\begin{aligned}
E[W(0)]=\widetilde{W}_{0}^{*'}(z)|_{z=1}
=\frac{n+1}{\mu_{b}}-1+\frac{\mu_{b}-\mu_{\nu}}{\theta\mu_{b}}\left[1-\left(\frac{\mu_{\nu}-\mu_{\nu}\theta}
{\theta+\mu_{\nu}-\theta\mu_{\nu}}\right)^{n+1}\right]+\frac{\theta+\mu_{b}-\theta\mu_{b}}{\mu_{b}(\theta+\mu_{\nu}-\theta\mu_{\nu})}.
\end{aligned}
\end{equation*}

\begin{figure}
\centering
\includegraphics[width=10cm,height=8cm]{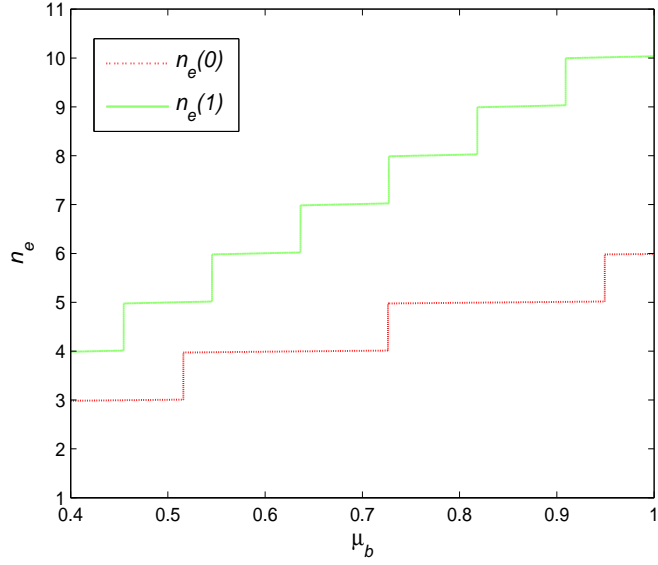}
\caption{Equilibrium thresholds for the observable systems with $R=10$, $C=1$, $\mu_{\nu}=0.4$, $\theta=0.2$}\label{Fig.1}
\end{figure}

Hence, the expected net benefit of the customer who joins the queueing system is given by
\begin{equation*}
\begin{aligned}
U(0)&=R-CE[W(0)]\\
&=R-C\left\{\frac{n+1}{\mu_{b}}-1+\frac{\mu_{b}-\mu_{\nu}}{\theta\mu_{b}}\left[1-\left(\frac{\mu_{\nu}-\mu_{\nu}\theta}
{\theta+\mu_{\nu}-\theta\mu_{\nu}}\right)^{n+1}\right]+\frac{\theta+\mu_{b}
-\theta\mu_{b}}{\mu_{b}(\theta+\mu_{\nu}-\theta\mu_{\nu})}\right\},
\end{aligned}
\end{equation*}
\begin{equation*}
U(1)=R-CE[W(1)]=R-C\left(\frac{n+1}{\mu_{b}}-1\right).
\end{equation*}
Solving $U(0)=0$ and $U(1)=0$, we obtain the positive and feasible
roots $(n_{e}^{*}(0),n_{e}^{*}(1))$. Because the uniqueness of the
roots in the observable case, we present a set of numerical
experiments to show the effect of the information level. We get that
the customers' equilibrium threshold strategy is
$(n_{e}(0),n_{e}(1))=(\lfloor n_{e}^{*}(0) \rfloor,\lfloor
n_{e}^{*}(1) \rfloor)$ in the observable case. From the
Figure~\ref{Fig.1} we can obtain that $n_{e}(0)$ and $n_{e}(1)$ both
increase with respect to $\mu_{b}$, however,
 $n_{e}(1)$ increases faster than $n_{e}(0)$.

\begin{figure}
\centering
\includegraphics[width=15cm,height=5.5cm]{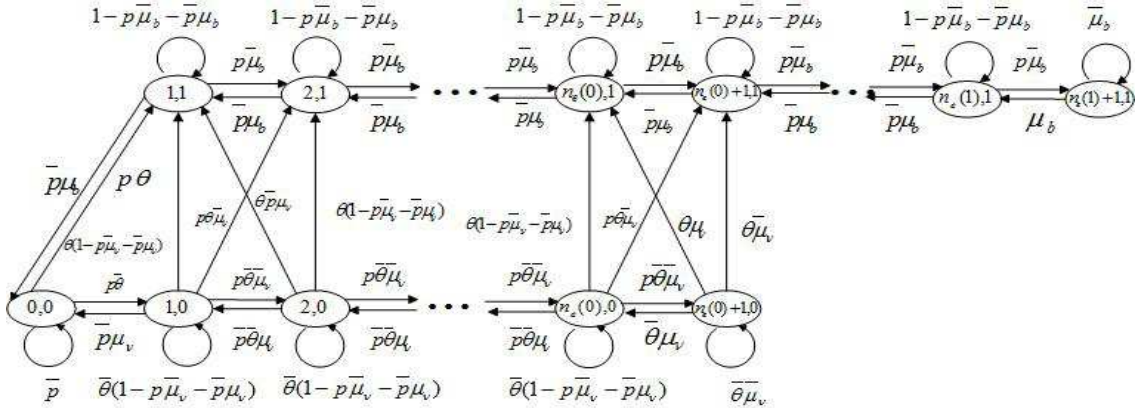}
\caption{Transition rate diagram for the observable queues}\label{Fig.2}
\end{figure}

Now we consider the stationary distribution in the observable case.
If all customers follow the threshold strategy $(n_{e}(0),n_{e}(1))$, the queueing system
conforms a Markov chain with state space
$\Omega_{ob}=\{(n,0)|0\leq n \leq n_{e}(0)+1\}\bigcup \{(n,1)|1\leq n \leq n_{e}(1)+1\}$.
The transition rate diagram is depicted in Figure~\ref{Fig.2}.  \\

Using the lexicographical ordering for the states, the transition probability matrix can be written as:
\[
        \widetilde{\mathbf{P}}=
         \bordermatrix{
         & \scriptstyle0 & \scriptstyle 1 & \cdots & \cdots & \scriptstyle n_{e}(0) & \scriptstyle n_{e}(0)+1 & \cdots & \cdots & \cdots & \cdots & \scriptstyle n_{e}(1)+1 \cr
         \scriptstyle0 & \mathbf{A}_{0} & \mathbf{C}_{0} & \cr
         \scriptstyle1 & \mathbf{B}_{0} & \mathbf{A}_{1} & \mathbf{C}_{1} & \cr
         \vdots &  & \mathbf{B}_{1} & \mathbf{A}_{1} & \mathbf{C}_{1} & \cr
         \vdots &  && \ddots & \ddots & \ddots & \cr
         \scriptstyle n_{e}(0) &  &&& \mathbf{B}_{1} & \mathbf{A}_{1} & \mathbf{C}_{1} & \cr
         \scriptstyle n_{e}(0)+1 &  &&&& \mathbf{B}_{2} & \mathbf{A}_{2} & \mathbf{C}_{2} & \cr
         \vdots &  &&&&& \mathbf{B}_{3} & \mathbf{A}_{3} & \mathbf{C}_{3} & \cr
         \vdots &  &&&&&& \mathbf{B}_{4} & \mathbf{A}_{3} & \mathbf{C}_{3} & \cr
         \vdots &  &&&&&&& \ddots & \ddots & \ddots & \cr
         \vdots &  &&&&&&&& \mathbf{B}_{4} & \mathbf{A}_{3} & \mathbf{C}_{3} \cr
         \scriptstyle n_{e}(1)+1 & &&&&&&&&& \mathbf{B}_{5} & \mathbf{A}_{4}
         },
         \]
where\\
$\mathbf{A}_{0}=\bar{p}$, $\mathbf{C}_{0}=[p\bar{\theta},p\theta]$,
$\mathbf{B}_{0}=\begin{bmatrix}
\bar{p}\mu_{\nu} \\
\bar{p}\mu_{b}
\end{bmatrix}$,
$
\mathbf{B}_{1}=\left[\begin{array}{cc}
       \bar{p}\bar{\theta}\mu_{\nu} & \bar{p}\theta\mu_{\nu} \\
       0 & \bar{p}\mu_{b}
       \end{array} \right],
$\\ $ \mathbf{A}_{1}=\left[\begin{array}{cc}
       \bar{\theta}(1-p\bar{\mu}_{\nu}-\bar{p}\mu_{\nu}) & \theta(1-p\bar{\mu}_{\nu}-\bar{p}\mu_{\nu}) \\
       0 & 1-p\bar{\mu}_{b}-\bar{p}\mu_{b}
       \end{array} \right], $
$
\mathbf{C}_{1}=\left[\begin{array}{cc}
       p\bar{\theta}\bar{\mu}_{\nu} & p\theta\bar{\mu}_{\nu} \\
       0 & p\bar{\mu}_{b}
       \end{array} \right], $\\
$
\mathbf{B}_{2}=\left[\begin{array}{cc}
       \bar{\theta}\mu_{\nu} & \theta\mu_{\nu} \\
       0 & \bar{p}\mu_{b}
       \end{array} \right], $
$
\mathbf{A}_{2}=\left[\begin{array}{cc}
       \bar{\theta}\bar{\mu}_{\nu} & \theta\bar{\mu}_{\nu} \\
       0 & 1-p\bar{\mu}_{b}-\bar{p}\mu_{b}
       \end{array} \right],
$
$\mathbf{C}_{2}=\begin{bmatrix}
0 \\
p\bar{\mu}_{b}
\end{bmatrix}$,\\
$\mathbf{B}_{3}=[0,\bar{p}\mu_{b}]$,
$\mathbf{A}_{3}=1-p\bar{\mu}_{b}-\bar{p}\mu_{b}$,
$\mathbf{C}_{3}=p\bar{\mu}_{b}$, $\mathbf{B}_{4}=\bar{p}\mu_{b}$,
$\mathbf{B}_{5}=\mu_{b}$, $\mathbf{A}_{4}=\bar{\mu}_{b}$.

Let $(L^{+},J)$ be the stationary limit of ${(L_{n}^{+},J_{n})}$ and the stationary distribution in the observable case is denoted as
\[
\pi_{nj}^{+}=P\{L^{+}=n,J=j\},\ \ (n,j)\in\Omega_{ob}. \]
\[
     \pi_{n}=\left\{\begin{array}{ll}
                    \pi_{00}^{+}, ~n=0,\\
                    (\pi_{n0}^{+},\pi_{n1}^{+}), ~1\leq n \leq n_{e}(0)+1,\\
                    \pi_{n1}^{+}, ~n_{e}(0)+2\leq n \leq n_{e}(1)+1.
                  \end{array}
\right.
\]
\[
\bm{\pi}=(\pi_{0},\pi_{1},...,\pi_{n_{e}(1)+1}).
\]
So the stationary transition probability equations can be written as
\begin{align}
&\pi_{00}^{+}=\bar{p}\pi_{00}^{+}+\bar{p}\mu_{\nu}\pi_{10}^{+}+\bar{p}\mu_{b}\pi_{11}^{+},\label{3.3}\\
&\pi_{10}^{+}=p\bar{\theta}\pi_{00}^{+}+\bar{\theta}(1-p\bar{\mu}_{\nu}-\bar{p}\mu_{\nu})\pi_{10}^{+}+\bar{\theta}\bar{p}\mu_{\nu}\pi_{20}^{+},\label{3.4}\\
&\pi_{n0}^{+}=\bar{\theta}p\bar{\mu}_{\nu}\pi_{n-1,0}^{+}+\bar{\theta}(1-p\bar{\mu}_{\nu}-\bar{p}\mu_{\nu})\pi_{n0}^{+}+\bar{\theta}\bar{p}\mu_{\nu}\pi_{n+1, 0}^{+},~~n=2,...,n_{e}(0)-1,\label{3.5}\\
&\pi_{n_{e}(0),0}^{+}=\bar{\theta}p\bar{\mu}_{\nu}\pi_{n_{e}(0)-1,0}^{+}+\bar{\theta}(1-p\bar{\mu}_{\nu}-\bar{p}\mu_{\nu})\pi_{n_{e}(0), 0}^{+}+\bar{\theta}\mu_{\nu}\pi_{n_{e}(0)+1,0}^{+}, \label{3.6}\\
&\pi_{n_{e}(0)+1,0}^{+}=p\bar{\theta}\bar{\mu}_{\nu}\pi_{n_{e}(0),0}^{+}+\bar{\theta}\bar{\mu}_{\nu}\pi_{n_{e}(0)+1,0}^{+},\label{3.7}\\
&\pi_{11}^{+}=p\theta\pi_{00}^{+}+\theta(1-p\bar{\mu}_{\nu}-\bar{p}\mu_{\nu})\pi_{10}^{+}+(1-p\bar{\mu}_{b}-\bar{p}\mu_{b})\pi_{11}^{+}+\theta\bar{p}\mu_{\nu}\pi_{20}^{+}+\bar{p}\mu_{b}\pi_{21}^{+},\label{3.8}
\end{align}
\begin{align}
&\pi_{n1}^{+}=\theta p\bar{\mu}_{\nu}\pi_{n-1,0}^{+}+p\bar{\mu}_{b}\pi_{n-1,1}^{+}+\theta(1-p\bar{\mu}_{\nu}-\bar{p}\mu_{\nu})\pi_{n0}^{+}
+(1-p\bar{\mu}_{b}-\bar{p}\mu_{b})\pi_{n1}^{+}+\theta\bar{p}\mu_{\nu}\pi_{n+1,0}^{+} \nonumber\\
&~~~~~~~~+\bar{p}\mu_{b}\pi_{n+1,1}^{+},~~n=2,...,n_{e}(0)-1,\label{3.9}\\
&\pi_{n_{e}(0),1}^{+}=\theta p\bar{\mu}_{\nu}\pi_{n_{e}(0)-1,0}^{+}+p\bar{\mu}_{b}\pi_{n_{e}(0)-1, 1}^{+}+\theta(1-p\bar{\mu}_{\nu}-\bar{p}\mu_{\nu})\pi_{n_{e}(0),0}^{+}+(1-p\bar{\mu}_{b}-\bar{p}\mu_{b})\pi_{n_{e}(0),1}^{+}  \nonumber \\
&~~~~~~~~~~~~+\theta\mu_{\nu}\pi_{n_{e}(0)+1,0}^{+}+\bar{p}\mu_{b}\pi_{n_{e}(0)+1,1}^{+},\label{3.10}\\
&\pi_{n_{e}(0)+1,1}^{+}=p\bar{\mu}_{b}\pi_{n_{e}(0),1}^{+}+\bar{p}\mu_{b}\pi_{n_{e}(0)+2,1}^{+}+(1-p\bar{\mu}_{b}-\bar{p}\mu_{b})\pi_{n_{e}(0)+1,1}^{+}
+\theta\bar{\mu}_{\nu}\pi_{n_{e}(0)+1,0}^{+}  \nonumber \\
&~~~~~~~~~~~~~~~+p\theta\bar{\mu}_{\nu}\pi_{n_{e}(0),0}^{+},\label{3.11}\\
&\pi_{n1}^{+}=p\bar{\mu}_{b}\pi_{n-1,1}^{+}+(1-p\bar{\mu}_{b}-\bar{p}\mu_{b})\pi_{n 1}^{+}+\bar{p}\mu_{b}\pi_{n+1, 1}^{+},~~n=n_{e}(0)+2,...,n_{e}(1)-1,\label{3.12}\\
&\pi_{n_{e}(1),1}^{+}=p\bar{\mu}_{b}\pi_{n_{e}(1)-1,1}^{+}+(1-p\bar{\mu}_{b}-\bar{p}\mu_{b})\pi_{n_{e}(1),1}^{+}+\mu_{b}\pi_{n_{e}(1)+1,1}^{+},\label{3.13}\\
&\pi_{n_{e}(1)+1,1}^{+}=p\bar{\mu}_{b}\pi_{n_{e}(1),1}^{+}+\bar{\mu}_{b}\pi_{n_{e}(1)+1,1}^{+}.\label{3.14}
\end{align}
Define $\alpha=\frac{p\bar{\mu}_{b}}{\bar{p}\mu_{b}}<1$.

Taking into account (\ref{3.3})--(\ref{3.7}), we first consider the probabilities $\{\pi_{n 0}^{+}|0\leq n \leq n_{e}(0)+1\}$. From (\ref{3.3}) we can get that
\begin{align}
\pi_{00}^{+}=\frac{\bar{p}\mu_{\nu}}{p}\pi_{10}^{+}+\frac{\bar{p}\mu_{b}}{p}\pi_{11}^{+}.\label{3.15}
\end{align}
Substituting (\ref{3.15}) into (\ref{3.4}), we obtain
\begin{align}
(\theta+\bar{\theta}p\bar{\mu}_{\nu})\pi_{10}^{+}=\bar{p}\bar{\theta}\mu_{b}\pi_{11}^{+}+\bar{\theta}\bar{p}\mu_{\nu}\pi_{20}^{+}.\label{3.16}
\end{align}
From (\ref{3.7}) we get
\begin{align}
\pi_{n_{e}(0)+1,0}^{+}=\frac{p\bar{\theta}\bar{\mu}_{\nu}}{1-\bar{\theta}\bar{\mu}_{\nu}}\pi_{n_{e}(0),0}^{+}.\label{3.17}
\end{align}
Substituting (\ref{3.17}) into (\ref{3.6}), we have
\begin{align}
(\theta+\bar{\theta}p\bar{\mu}_{\nu}+\bar{\theta}\bar{p}\mu_{\nu}-\theta\bar{\theta}\bar{\mu}_{\nu}-{\bar{\theta}}^{2}{\bar{\mu}_{\nu}}^{2}p-{\bar{\theta}}^{2}\bar{\mu}_{\nu}\mu_{\nu})\pi_{n_{e}(0), 0}^{+}-p\bar{\theta}\bar{\mu}_{\nu}(1-\bar{\theta}\bar{\mu}_{\nu})\pi_{n_{e}(0)-1,0}^{+}=0.\label{3.18}
\end{align}
From (\ref{3.5}) we find that $\{\pi_{n 0}^{+}|1\leq n \leq n_{e}(0)\}$ are the solutions of the following homogeneous linear difference equation with constant coefficients:
\begin{align}
\bar{\theta}\bar{p}\mu_{\nu}x_{n+1}-(\theta+\bar{\theta}p\bar{\mu}_{\nu}+\bar{\theta}\bar{p}\mu_{\nu})x_{n}+\bar{\theta}p\bar{\mu}_{\nu}x_{n-1}=0,~~n=2,3,...,n_{e}(0)-1.\label{3.19}
\end{align}
We investigate the corresponding characteristic equation
\begin{align}
\bar{\theta}\bar{p}\mu_{\nu}x^{2}-(\theta+\bar{\theta}p\bar{\mu}_{\nu}+\bar{\theta}\bar{p}\mu_{\nu})x+\bar{\theta}p\bar{\mu}_{\nu}=0, \nonumber
\end{align}
which has two roots
\begin{align}
x_{1,2}^{*}=\frac{(\theta+\bar{\theta}p\bar{\mu}_{\nu}+\bar{\theta}\bar{p}\mu_{\nu})\pm\sqrt{(\theta+\bar{\theta}p\bar{\mu}_{\nu}+\bar{\theta}\bar{p}\mu_{\nu})^{2}-4{\bar{\theta}}^{2}p\bar{p}\mu_{\nu}\bar{\mu}_{\nu}}}{2\bar{\theta}\bar{p}\mu_{\nu}}.\label{3.20}
\end{align}
The general solution of (\ref{3.19}), denoted by $x_{n}^{hom}$, is $x_{n}^{hom}=\widetilde{A}_{1}x_{1}^{*n}+\widetilde{B}_{1}x_{2}^{*n}~(x_{1}^{*}\neq x_{2}^{*})$, where $\widetilde{A}_{1}$, $\widetilde{B}_{1}$ are the coefficients to be determined. From (\ref{3.16}) and (\ref{3.18}) we get the equations about $\widetilde{A}_{1}$ and $\widetilde{B}_{1}$ as
\begin{equation}\label{3.21}
     \left\{
     \begin{aligned}
                    &\left(\theta x_{1}^{*}+\bar{\theta}p\bar{\mu}_{\nu}x_{1}^{*}-
                    \bar{\theta}\bar{p}\mu_{\nu}x_{1}^{*2}\right)\widetilde{A}_{1}+
                    \left(\theta x_{2}^{*}+\bar{\theta}p\bar{\mu}_{\nu}x_{2}^{*}-\bar{\theta}
                    \bar{p}\mu_{\nu}x_{2}^{*2}\right)\widetilde{B}_{1}=\bar{p}\bar{\theta}\mu_{b}\pi_{11}^{+},\\
                    &\left[\left(\theta+\bar{\theta}p\bar{\mu}_{\nu}+\bar{\theta}\bar{p}\mu_{\nu}-
                    \theta\bar{\theta}\bar{\mu}_{\nu}-{\bar{\theta}}^{2}{\bar{\mu}_{\nu}}^{2}p
                    -{\bar{\theta}}^{2}\bar{\mu}_{\nu}\mu_{\nu}\right)x_{1}^{n_{e}(0)}-
                    p\bar{\theta}\bar{\mu}_{\nu}\left(1-\bar{\theta}\bar{\mu}_{\nu}\right)x_{1}^{n_{e}(0)-1}\right]
                    \widetilde{A}_{1}\\ &+\left[\left(\theta+\bar{\theta}p\bar{\mu}_{\nu}+\bar{\theta}\bar{p}
                    \mu_{\nu}-\theta\bar{\theta}\bar{\mu}_{\nu}-{\bar{\theta}}^{2}{\bar{\mu}_{\nu}}^{2}p
                    -{\bar{\theta}}^{2}\bar{\mu}_{\nu}\mu_{\nu}\right)x_{2}^{n_{e}(0)}-p\bar{\theta}
                    \bar{\mu}_{\nu}\left(1-\bar{\theta}\bar{\mu}_{\nu}\right)x_{2}^{n_{e}(0)-1}\right]\widetilde{B}_{1}=0.
                    \end{aligned}
\right.
\end{equation}
Solving (\ref{3.21}), we obtain
\begin{equation}\label{3.22}
     \left\{
     \begin{aligned}
                   \widetilde{A}_{1}=&H\bar{p}\bar{\theta}\mu_{b}\pi_{11}^{+}\\
                   &\cdot\left[p\bar{\theta}\bar{\mu}_{\nu}
                   \left(1-\bar{\theta}\bar{\mu}_{\nu}\right)x_{2}^{*n_{e}(0)-1}-
                   \left(\theta+\bar{\theta}p\bar{\mu}_{\nu}+\bar{\theta}\bar{p}\mu_{\nu}
                   -\theta\bar{\theta}\bar{\mu}_{\nu}-{\bar{\theta}}^{2}{\bar{\mu}_{\nu}}^{2}p
                   -{\bar{\theta}}^{2}\bar{\mu}_{\nu}\mu_{\nu}\right)x_{2}^{*n_{e}(0)}\right],\\
                   \widetilde{B}_{1}=&H\bar{p}\bar{\theta}\mu_{b}\pi_{11}^{+}\\
                   &\cdot\left[\left(\theta+\bar{\theta}p
                   \bar{\mu}_{\nu}+\bar{\theta}\bar{p}\mu_{\nu}-\theta\bar{\theta}\bar{\mu}_{\nu}
                   -{\bar{\theta}}^{2}{\bar{\mu}_{\nu}}^{2}p-{\bar{\theta}}^{2}\bar{\mu}_{\nu}\mu_{\nu}\right)
                   x_{1}^{*n_{e}(0)}-p\bar{\theta}\bar{\mu}_{\nu}\left(1-\bar{\theta}\bar{\mu}_{\nu}\right)
                   x_{1}^{*n_{e}(0)-1}\right],
     \end{aligned}
\right.
\end{equation}
where
\begin{align*}
H=&\left\{\left(\theta
x_{2}^{*}+\bar{\theta}p\bar{\mu}_{\nu}x_{2}^{*}-\bar{\theta}\bar{p}\mu_{\nu}x_{2}^{*2}\right)
\left[\left(\theta+\bar{\theta}p\bar{\mu}_{\nu}+\bar{\theta}\bar{p}\mu_{\nu}-\theta\bar{\theta}\bar{\mu}_{\nu}
-{\bar{\theta}}^{2}{\bar{\mu}_{\nu}}^{2}p-{\bar{\theta}}^{2}\bar{\mu}_{\nu}\mu_{\nu}\right)x_{1}^{*n_{e}(0)}\right.\right.\\
&\left.-p\bar{\theta}\bar{\mu}_{\nu}\left(1-\bar{\theta}\bar{\mu}_{\nu}\right)
x_{1}^{*n_{e}(0)-1}\right]-\left(\theta
x_{1}^{*}+\bar{\theta}p\bar{\mu}_{\nu}x_{1}^{*}-\bar{\theta}\bar{p}\mu_{\nu}x_{1}^{*2}\right)\\
&\left.\cdot
\left[\left(\theta+\bar{\theta}p\bar{\mu}_{\nu}+\bar{\theta}\bar{p}\mu_{\nu}-\theta\bar{\theta}\bar{\mu}_{\nu}
-{\bar{\theta}}^{2}{\bar{\mu}_{\nu}}^{2}p-{\bar{\theta}}^{2}\bar{\mu}_{\nu}\mu_{\nu}\right)x_{2}^{*n_{e}(0)}
-p\bar{\theta}\bar{\mu}_{\nu}(1-\bar{\theta}\bar{\mu}_{\nu})
x_{2}^{*n_{e}(0)-1}\right]\right\}^{-1}.
\end{align*} Thus,
\begin{equation*}
     \left\{
     \begin{aligned}
                    &\pi_{00}^{+}=\frac{\bar{p}\mu_{\nu}}{p}(\widetilde{A}_{1}x_{1}^{*}+\widetilde{B}_{1}x_{2}^{*})+\frac{\bar{p}\mu_{b}}{p}\pi_{11}^{+},\\
                    &\pi_{n0}^{+}=\widetilde{A}_{1}x_{1}^{*n}+\widetilde{B}_{1}x_{2}^{*n},~~n=1,2,...,n_{e}(0), \\
                    &\pi_{n_{e}(0)+1, 0}^{+}=\frac{p\bar{\theta}\bar{\mu}_{\nu}}{1-\bar{\theta}\bar{\mu}_{\nu}}(\widetilde{A}_{1}x_{1}^{*n_{e}(0)}+\widetilde{B}_{1}x_{2}^{*n_{e}(0)}).
     \end{aligned}
\right.
\end{equation*}

Next, we discuss the probabilities
$\{\pi_{n 1}^{+}|1\leq n \leq n_{e}(0)\}$. From (\ref{3.9}),
$\{\pi_{n 1}^{+}|1\leq n \leq n_{e}(0)\}$ are the solutions of the following nonhomogeneous linear difference equation:
\begin{align}
&\bar{p}\mu_{b}\pi_{n+1,1}^{+}-(p\bar{\mu}_{b}+\bar{p}\mu_{b})\pi_{n 1}^{+}+p\bar{\mu}_{b}\pi_{n-1,1}^{+}
=-\theta p\bar{\mu}_{\nu}(\widetilde{A}_{1}x_{1}^{*n-1}+\widetilde{B}_{1}x_{2}^{*n-1})\nonumber\\
&-\theta(1-p\bar{\mu}_{\nu}-\bar{p}\mu_{\nu})(\widetilde{A}_{1}x_{1}^{*n}+\widetilde{B}_{1}x_{2}^{*n})-\theta\bar{p}\mu_{\nu}(\widetilde{A}_{1}x_{1}^{*n+1}+\widetilde{B}_{1}x_{2}^{*n+1}), \nonumber \\
&n=2,3,...,n_{e}(0)-1.\label{3.23}
\end{align}
The general solution of the homogeneous version of (\ref{3.23}) is $x_{n}^{hom}=\widetilde{A}_{2}1^{n}+\widetilde{B}_{2}\alpha^{n}$. So the general solution of (\ref{3.23}), denoted by $x_{n}^{gen}$, is given as $x_{n}^{gen}=x_{n}^{hom}+x_{n}^{spec}$, where $x_{n}^{spec}$ is a specific solution of (\ref{3.23}). Because the nonhomogeneous part of (\ref{3.23}) is geometric with parameter $x_{1}^{*}$ and $x_{2}^{*}$, we consider a specific solution of the form $x_{n}^{spec}=\widetilde{C}_{1}x_{1}^{*n}+\widetilde{D}_{1}x_{2}^{*n}$. Substituting $x_{n}^{spec}=\widetilde{C}_{1}x_{1}^{*n}+\widetilde{D}_{1}x_{2}^{*n}$ into (\ref{3.23}), we get
\begin{equation}\label{3.24}
     \left\{
     \begin{aligned}
                    &\widetilde{C}_{1}=\frac{\theta \widetilde{A}_{1}[(x_{1}^{*}-1)(\bar{p}\mu_{\nu}x_{1}^{*}-p\bar{\mu}_{\nu})+x_{1}^{*}]}{(x_{1}^{*}-1)(p\bar{\mu}_{b}-\bar{p}\mu_{b}x_{1}^{*})},\\
                    &\widetilde{D}_{1}=\frac{\theta \widetilde{B}_{1}[(x_{2}^{*}-1)(\bar{p}\mu_{\nu}x_{2}^{*}-p\bar{\mu}_{\nu})+x_{2}^{*}]}{(1-x_{2}^{*})(\bar{p}\mu_{b}x_{2}^{*}-p\bar{\mu}_{b})}.
     \end{aligned}
\right.
\end{equation}
Hence, the general solution of (\ref{3.23}) is given as
\begin{align}
x_{n}^{gen}=\widetilde{A}_{2}1^{n}+\widetilde{B}_{2}\alpha^{n}+\widetilde{C}_{1}x_{1}^{*n}+\widetilde{D}_{1}x_{2}^{*n},~~n=1,2,...,n_{e}(0), \nonumber
\end{align}
where $\widetilde{A}_{2}$, $\widetilde{B}_{2}$ are to be determined. Taking into account (\ref{3.8}), we get that
\begin{equation}\label{3.25}
     \left\{
     \begin{aligned}
                    &\widetilde{A}_{2}+\alpha \widetilde{B}_{2}=\pi_{11}^{+}-\widetilde{C}_{1}x_{1}^{*}-\widetilde{D}_{1}x_{2}^{*},\\
                    &\bar{p}\mu_{b}(\widetilde{A}_{2}+\widetilde{B}_{2}\alpha^{2})=(p\bar{\mu}_{b}+\bar{p}\mu_{b}\bar{\theta})\pi_{11}^{+}-\theta(1-p\bar{\mu}_{\nu})(\widetilde{A}_{1}x_{1}^{*}+\widetilde{B}_{1}x_{2}^{*})\\
                    &~~~~~~~~~~~~~~~~~~~~~~~~-\theta\bar{p}\mu_{\nu}(\widetilde{A}_{1}x_{1}^{*2}+\widetilde{B}_{1}x_{2}^{*2})-\bar{p}\mu_{b}(\widetilde{C}_{1}x_{1}^{*2}+\widetilde{D}_{1}x_{2}^{*2}).
     \end{aligned}
\right.
\end{equation}
Solving (\ref{3.25}), we have
\begin{equation}\label{3.26}
      \left\{
      \begin{aligned}
                     &\widetilde{A}_{2}=\textstyle\frac{\bar{p}\mu_{b}\bar{\theta}\pi_{11}^{+}-\bar{p}\mu_{b}(\widetilde{C}_{1}x_{1}^{*2}+\widetilde{D}_{1}x_{2}^{*2}-\alpha \widetilde{C}_{1}x_{1}^{*}-\alpha
\widetilde{D}_{1}x_{2}^{*})-\theta\bar{p}\mu_{\nu}(\widetilde{A}_{1}x_{1}^{*2}+\widetilde{B}_{1}x_{2}^{*2})-\theta(1-p\bar{\mu}_{\nu})(\widetilde{A}_{1}x_{1}^{*}+\widetilde{B}_{1}x_{2}^{*})}{\bar{p}\mu_{b}(1-\alpha)},\\
                     &\widetilde{B}_{2}=\textstyle
\frac{\bar{p}\mu_{b}(\widetilde{C}_{1}x_{1}^{*2}+\widetilde{D}_{1}x_{2}^{*2}-\widetilde{C}_{1}x_{1}^{*}-\widetilde{D}_{1}x_{2}^{*})+\theta\bar{p}\mu_{\nu}(\widetilde{A}_{1}x_{1}^{*2}+\widetilde{B}_{1}x_{2}^{*2})
+\theta(1-p\bar{\mu}_{\nu})(\widetilde{A}_{1}x_{1}^{*}+\widetilde{B}_{1}x_{2}^{*})-(p\bar{\mu}_{b}-\bar{p}\mu_{b}\theta)\pi_{11}^{+}}{\bar{p}\mu_{b}\alpha(1-\alpha)}.
\end{aligned}
\right.
\end{equation}
Thus,
\begin{align}
\pi_{n1}^{+}=\widetilde{A}_{2}+\widetilde{B}_{2}\alpha^{n}+\widetilde{C}_{1}x_{1}^{*n}+\widetilde{D}_{1}x_{2}^{*n},~~n=1,2,...,n_{e}(0). \nonumber
\end{align}

Then we discuss the probabilities $\{\pi_{n 1}^{+}|n_{e}(0)+2\leq n \leq n_{e}(1)+1\}$. From (\ref{3.12}), $\{\pi_{n 1}^{+}|n_{e}(0)+2\leq n \leq n_{e}(1)\}$ are the solutions of the homogeneous version of (\ref{3.23}), i.e.,$x_{n}^{hom}=\widetilde{A}_{3}1^{n}+\widetilde{B}_{3}\alpha^{n}$, where $\widetilde{A}_{3}$, $\widetilde{B}_{3}$ are to be determined.\\
Substituting $x_{n}^{hom}=\widetilde{A}_{3}1^{n}+\widetilde{B}_{3}\alpha^{n}$ into (\ref{3.13}), we get that
\begin{equation}\label{3.27}
      \left\{
      \begin{aligned}
                    &\widetilde{A}_{3}=0,\\
                    &\widetilde{B}_{3}=\alpha^{-1-n_{e}(0)}\pi_{n_{e}(0)+1,1}^{+}.
      \end{aligned}
\right.
\end{equation}
From (\ref{3.14}), we obtain
\begin{equation*}
\pi_{n_{e}(1)+1,1}=\bar{p}\alpha^{n_{e}(1)-n_{e}(0)}\pi_{n_{e}(0)+1,1}^{+}.
\end{equation*}
Thus,
\begin{equation*}
     \left\{
     \begin{aligned}
                    &\pi_{n1}^{+}=\widetilde{B}_{3}\alpha^{n},~~n=n_{e}(0)+2,n_{e}(0)+3,...,n_{e}(1),\\
                    &\pi_{n_{e}(1)+1,1}^{+}=\bar{p}\alpha^{n_{e}(1)-n_{e}(0)}\pi_{n_{e}(0)+1,1}^{+}.
     \end{aligned}
\right.
\end{equation*}

Finally, we consider the probability $\pi_{n_{e}(0)+1,1}^{+}$.\\
From (\ref{3.11}) we get
\begin{equation*}
\pi_{n_{e}(0)+1,                                                                                                                                                                                1}^{+}=\alpha(\widetilde{A}_{2}+\widetilde{B}_{2}\alpha^{n_{e}(0)}+\widetilde{C}_{1}x_{1}^{*n_{e}(0)}+\widetilde{D}_{1}x_{2}^{*n_{e}(0)})+\frac{p\theta\bar{\mu}_{\nu}}{\bar{p}\mu_{b}(1-\bar{\theta}\bar{\mu}_{\nu})}(\widetilde{A}_{1}x_{1}^{*n_{e}(0)}+\widetilde{B}_{1}x_{2}^{*n_{e}(0)}).
\end{equation*}
Thus,
\begin{equation*}
\left\{
\begin{aligned}
&\pi_{n1}^{+}=\widetilde{A}_{2}+\widetilde{B}_{2}\alpha^{n}+\widetilde{C}_{1}x_{1}^{*n}+\widetilde{D}_{1}x_{2}^{*n},~~n=1,2,...,n_{e}(0),\\
&\pi_{n_{e}(0)+1, 1}^{+}=\alpha(\widetilde{A}_{2}+\widetilde{B}_{2}\alpha^{n_{e}(0)}+\widetilde{C}_{1}x_{1}^{*n_{e}(0)}+\widetilde{D}_{1}x_{2}^{*n_{e}(0)})+\frac{p\theta\bar{\mu}_{\nu}}{\bar{p}\mu_{b}(1-\bar{\theta}\bar{\mu}_{\nu})}(\widetilde{A}_{1}x_{1}^{*n_{e}(0)}+\widetilde{B}_{1}x_{2}^{*n_{e}(0)}),\\
&\pi_{n1}^{+}=\widetilde{B}_{3}\alpha^{n},~~n=n_{e}(0)+2,n_{e}(0)+3,...,n_{e}(1),\\
&\pi_{n_{e}(1)+1, 1}^{+}=\bar{p}\alpha^{n_{e}(1)-n_{e}(0)}\left[\alpha(\widetilde{A}_{2}+\widetilde{B}_{2}\alpha^{n_{e}(0)+1}+\widetilde{C}_{1}x_{1}^{*n_{e}(0)+1}+\widetilde{D}_{1}x_{2}^{*n_{e}(0)+1})\right.\\
&~~~~~~~~~~~~~~\left.+\frac{p\theta\bar{\mu}_{\nu}}{\bar{p}\mu_{b}(1-\bar{\theta}\bar{\mu}_{\nu})}(\widetilde{A}_{1}x_{1}^{*n_{e}(0)}+\widetilde{B}_{1}x_{2}^{*n_{e}(0)})\right].
\end{aligned}
\right.
\end{equation*}

In conclusion, we have expressed all stationary probabilities in terms of $\pi_{11}^{+}$. The remaining probability $\pi_{11}^{+}$ can be solved by the normalization equation $\sum_{n=0}^{n_{e}(0)+1}\pi_{n0}^{+}+\sum_{n=1}^{n_{e}(1)+1}\pi_{n1}^{+}=1$, so we obtain the following theorem which provides all stationary probabilities.

\begin{theorem} \label{the:3.1}
Consider an observable $Geo/Geo/1$ queue with multiple working vocations, in which all arriving customers follow the threshold policy $(n_{e}(0),n_{e}(1))$. If
$\alpha<1$, the stationary distribution $\{\pi_{nj}^{+}|(n,j)\in\Omega_{ob}\}$ of $(L^{+},J)$ is
\begin{align}
&\pi_{00}^{+}=\frac{\bar{p}\mu_{\nu}}{p}(\widetilde{A}_{1}x_{1}^{*}+\widetilde{B}_{1}x_{2}^{*})+\frac{\bar{p}\mu_{b}}{p}\pi_{11}^{+},\\
&\pi_{n0}^{+}=\widetilde{A}_{1}x_{1}^{*n}+\widetilde{B}_{1}x_{2}^{*n},~~n=1,2,...,n_{e}(0),\\
&\pi_{n_{e}(0)+1, 0}^{+}=\frac{p\bar{\theta}\bar{\mu}_{\nu}}{1-\bar{\theta}\bar{\mu}_{\nu}}(\widetilde{A}_{1}x_{1}^{*n_{e}(0)}+\widetilde{B}_{1}x_{2}^{*n_{e}(0)}),\\
&\pi_{n1}^{+}=\widetilde{A}_{2}+\widetilde{B}_{2}\alpha^{n}+\widetilde{C}_{1}x_{1}^{*n}+\widetilde{D}_{1}x_{2}^{*n},~~n=1,2,...,n_{e}(0),\\
&\pi_{n_{e}(0)+1, 1}^{+}=\alpha(\widetilde{A}_{2}+\widetilde{B}_{2}\alpha^{n_{e}(0)}+\widetilde{C}_{1}x_{1}^{*n_{e}(0)}+\widetilde{D}_{1}x_{2}^{*n_{e}(0)})+\frac{p\theta\bar{\mu}_{\nu}}{\bar{p}\mu_{b}(1-\bar{\theta}\bar{\mu}_{\nu})}(\widetilde{A}_{1}x_{1}^{*n_{e}(0)}+\widetilde{B}_{1}x_{2}^{*n_{e}(0)}),\\
&\pi_{n1}^{+}=\widetilde{B}_{3}\alpha^{n},~~n=n_{e}(0)+2,n_{e}(0)+3,...,n_{e}(1),\\
&\pi_{n_{e}(1)+1, 1}^{+}=\bar{p}\alpha^{n_{e}(1)-n_{e}(0)}\left[\alpha(\widetilde{A}_{2}+\widetilde{B}_{2}\alpha^{n_{e}(0)+1}+\widetilde{C}_{1}x_{1}^{*n_{e}(0)+1}+\widetilde{D}_{1}x_{2}^{*n_{e}(0)+1})\right.\\ \nonumber
&~~~~~~~~~~~~~~\left.+\frac{p\theta\bar{\mu}_{\nu}}{\bar{p}\mu_{b}(1-\bar{\theta}\bar{\mu}_{\nu})}(\widetilde{A}_{1}x_{1}^{*n_{e}(0)}+\widetilde{B}_{1}x_{2}^{*n_{e}(0)})\right].
\end{align}
where $x_{1}^{*}$, $x_{2}^{*}$, $\widetilde{A}_{1}$, $\widetilde{B}_{1}$, $\widetilde{C}_{1}$, $\widetilde{D}_{1}$, $\widetilde{A}_{2}$, $\widetilde{B}_{2}$, $\widetilde{B}_{3}$ are given by (\ref{3.20}), (\ref{3.22}), (\ref{3.24}), (\ref{3.26}), (\ref{3.27}) respectively, and $\pi_{11}^{+}$ can be derived by the normalization equation $\sum_{n=0}^{n_{e}(0)+1}\pi_{n0}^{+}+\sum_{n=1}^{n_{e}(1)+1}\pi_{n1}^{+}=1$.
\end{theorem}

According to Theorem~\ref{the:3.1}, the probability of balking is equal to $\pi_{n_{e}(0)+1,0}^{+}+\pi_{n_{e}(1)+1,1}^{+}$. So the social benefit per time unit for the
threshold policy $(n_{e}(0),n_{e}(1))$ can be calculated as
\begin{equation*}
U_{s}(n_{e}(0),n_{e}(1))=pR(1-\pi_{n_{e}(0)+1,0}^{+}-\pi_{n_{e}(1)+1, 1}^{+})-C\left(\sum_{n=1}^{n_{e}(0)+1}n\pi_{n0}^{+}+\sum_{n=1}^{n_{e}(1)+1}n\pi_{n1}^{+}\right).
\end{equation*}

\begin{figure}
\centering
\includegraphics[width=10cm,height=8cm]{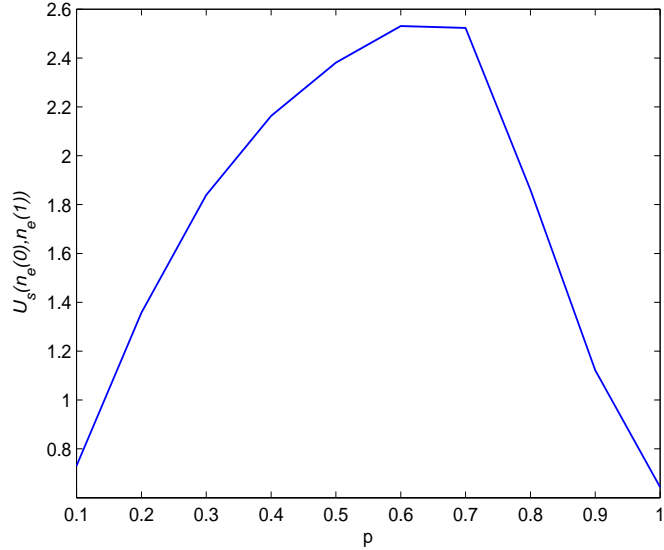}
\caption{Equilibrium social benefit for the observable systems with $R=10$, $C=1$, $\theta=0.05$, $\mu_{\nu}=0.4$, $\mu_{b}=0.8$}\label{Fig.3}
\end{figure}

When all customers follow the above equilibrium threshold strategy $(n_{e}(0),n_{e}(1))$, the social benefit per time unit in equilibrium can be represented
as $U_{s}(n_{e}(0),n_{e}(1))$. Fig.~\ref{Fig.3} is concerned with the social benefit under the equilibrium threshold strategy. We can observe that
$U_{s}(n_{e}(0),n_{e}(1))$ first increases, then decreases with respect to $p$.

\begin{figure}
\centering
\includegraphics[width=10cm,height=8cm]{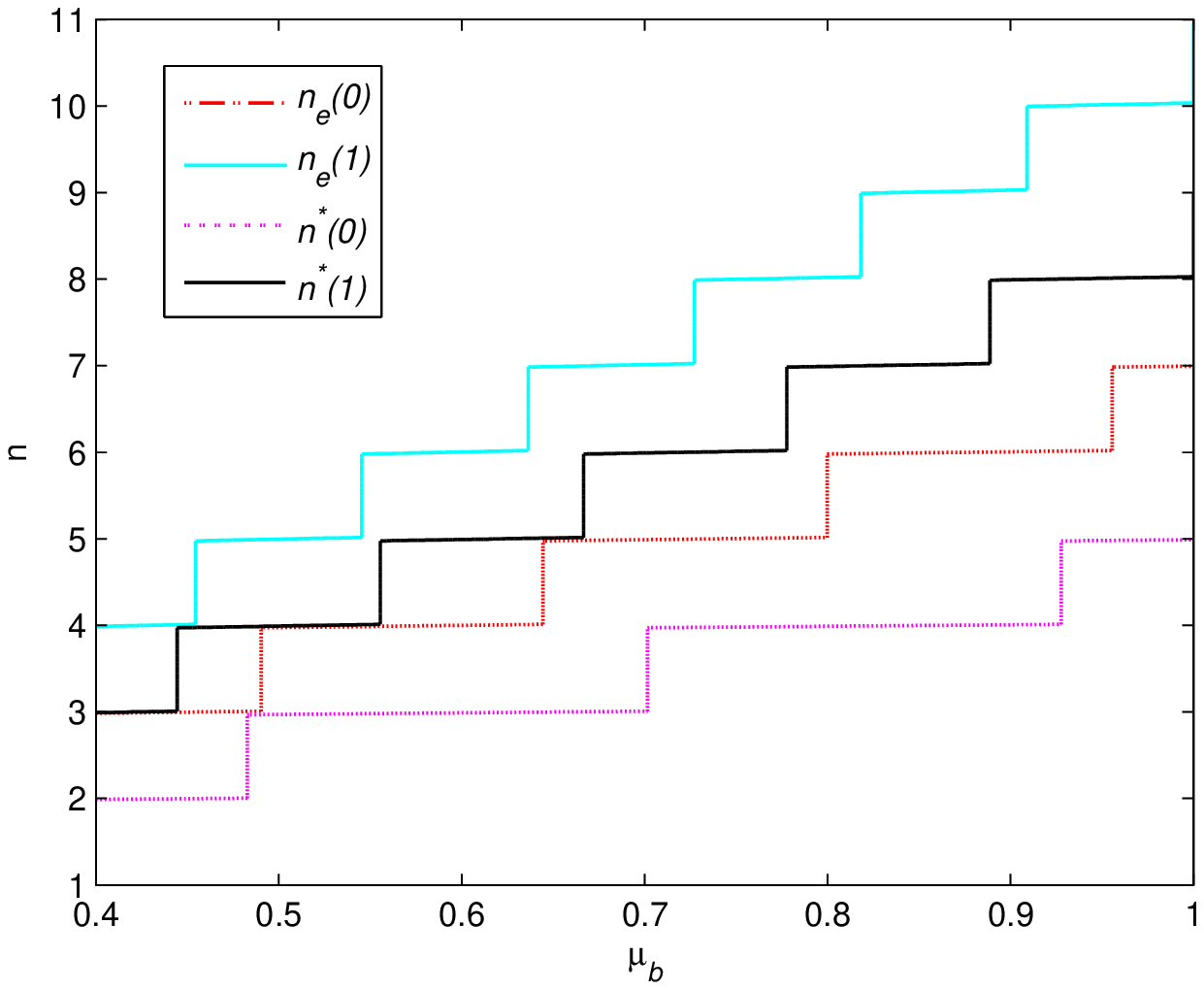}
\caption{Comparisons of equilibrium and socially optimal threshold strategies for the observable systems with $R=10$, $C=1$, $\theta=0.3$, $\mu_{\nu}=0.4$, $p=0.5$}\label{Fig.4}
\end{figure}

Next, from the view of social optimization, denote the socially
optimal threshold strategy as $(n^{*}(0),n^{*}(1))$, which can be
obtained by solving the unconstrained integer programming $\max
U_{s}(n_{e}(0),n_{e}(1))$. Figure~\ref{Fig.4} compares the
equilibrium threshold strategy $(n_{e}(0),n_{e}(1))$ and the
socially optimal threshold strategy $(n^{*}(0),n^{*}(1))$ for the
observable systems. We get that $n_{e}(0)>n^{*}(0)$ and
$n_{e}(1)>n^{*}(1)$, which shows that the individual optimization
deviates from the social optimization for the observable systems.

\section{The partially observable queue}
\label{4}

\begin{figure}
\centering
\includegraphics[width=15cm,height=5.5cm]{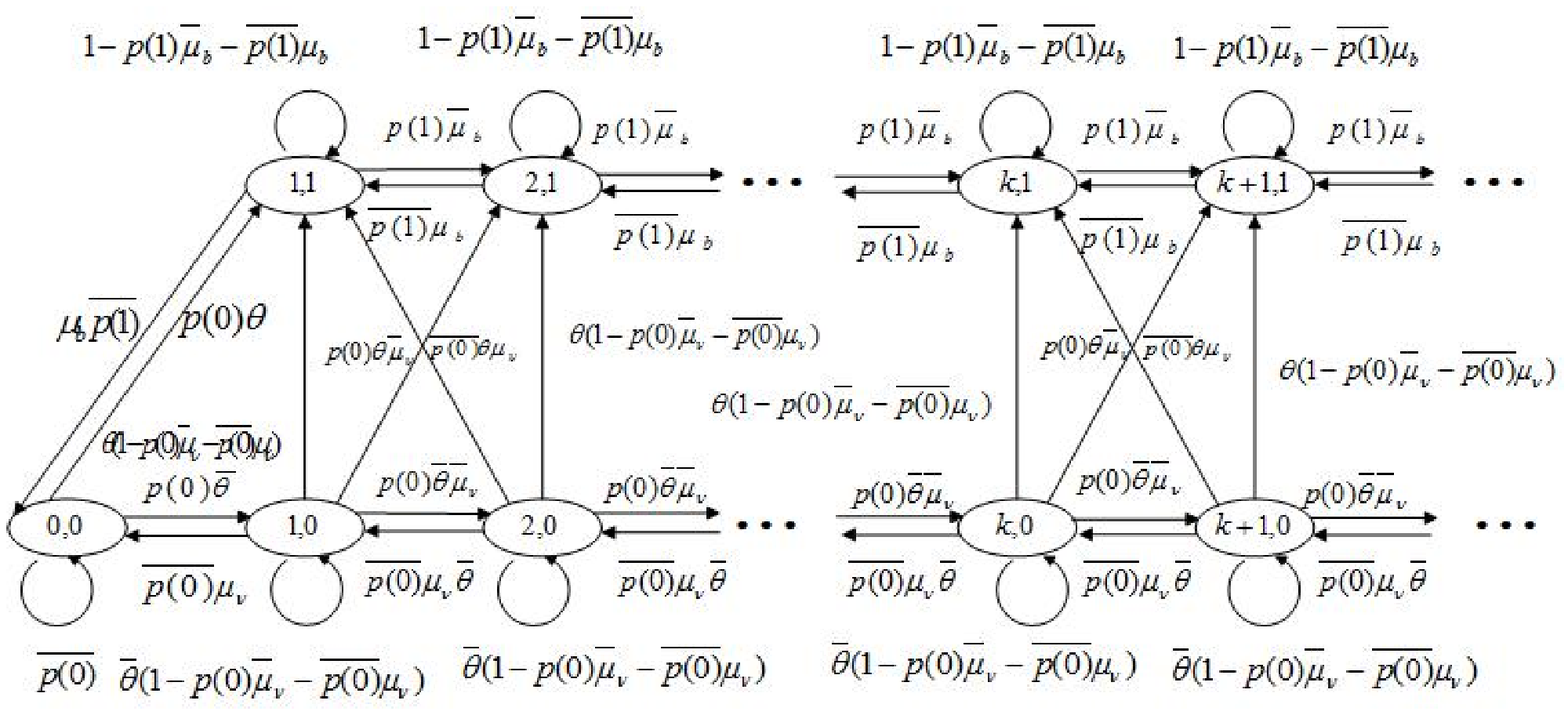}
\caption{Transition rate diagram for the partially observable queues}\label{Fig.5}
\end{figure}

In this section, we study the partially observable case in which the
arriving customers only observe the state $J_{n}$ of the server.
Suppose that the customers' decision that whether to join or to balk
upon their arrival can be indicated by a couple of joining
probabilities $(q(0),q(1))~(0\leq q(j)\leq1,~ j=0,1)$, i.e., the
arrival rate is $p(j)=pq(j)$ when the server is in state $j$, and
their equilibrium mixed strategy is represented by
$(q_{e}(0),q_{e}(1))$. The process
$\{(L_{n}^{+},J_{n}),n\geqslant0\}$ conforms a Markov chain with
state space
$\Omega_{po}=\{(0,0)\}\bigcup\{(k,j):k\geqslant1,~j=0,1\}$. The
transition rate diagram is described in Figure~\ref{Fig.5}.

Using the lexicographical ordering for the states, the transition probability matrix can be written as:
\[
        \mathbf{P}=
         \left[\begin{array}{cccccc}
         \mathbf{B}_{0}^{'} & \mathbf{C}_{0}^{'} & \\
         \mathbf{B}_{1}^{'} & \mathbf{A}_{1}^{'} & \mathbf{C}_{1}^{'} & \\
         & \mathbf{B}_{2}^{'} & \mathbf{A}_{1}^{'} & \mathbf{C}_{1}^{'} & \\
         && \mathbf{B}_{2}^{'} & \mathbf{A}_{1}^{'} & \mathbf{C}_{1}^{'} & \\
         &&& \ddots & \ddots & \ddots
        \end{array} \right],
         \]
where\\
$\mathbf{B}_{0}^{'}=\overline{p(0)}$, $\mathbf{C}_{0}^{'}=[p(0)\bar{\theta},p(0)\theta]$,
$\mathbf{B}_{1}^{'}=\begin{bmatrix}
\overline{p(0)}\mu_{\nu} \\
\overline{p(1)}\mu_{b}
\end{bmatrix}$,\\
$
\mathbf{A}_{1}^{'}=\left[\begin{array}{cc}
       \bar{\theta}(1-p(0)\bar{\mu}_{\nu}-\overline{p(0)}\mu_{\nu}) & \theta(1-p(0)\bar{\mu}_{\nu}-\overline{p(0)}\mu_{\nu}) \\
       0 & 1-p(1)\bar{\mu}_{b}-\overline{p(1)}\mu_{b}
       \end{array} \right],
$\\
$
\mathbf{C}_{1}^{'}=\left[\begin{array}{cc}
       p(0)\bar{\theta}\bar{\mu}_{\nu} & p(0)\theta\bar{\mu}_{\nu} \\
       0 & p(1)\bar{\mu}_{b}
       \end{array} \right],
$
$
\mathbf{B}_{2}^{'}=\left[\begin{array}{cc}
       \overline{p(0)}\bar{\theta}\mu_{\nu} & \overline{p(0)}\theta\mu_{\nu} \\
       0 & \overline{p(1)}\mu_{b}
       \end{array} \right].
$

Due to the block tridiagonal structure of the transition probability matrix, ${(L_{n}^{+},J_{n})}$ is a quasi-birth-and-death chain. Setting $\widetilde{\alpha}=\frac{p(1)\bar{\mu}_{b}}{\mu_{b}\overline{p(1)}}<1$, let $(L^{+},J)$ be the stationary limit of ${(L_{n}^{+},J_{n})}$ and its distribution is denoted as
\[
\pi_{kj}^{+'}=P\{L^{+}=k,J=j\},\ \ (k,j)\in\Omega_{po},
\]
\[
\pi_{00}^{+'}=\pi_{0}^{'}, \ \ (\pi_{k0}^{+'},\pi_{k1}^{+'})=\pi_{k}^{'}.
\]

Using the matrix-geometric solution method, we get that the rate matric $\mathbf{R}$ is the minimal non-negative solution of the matrix quadratic equation:
\begin{equation}\label{4.1}
\mathbf{R}=\mathbf{R}^{2}\mathbf{B}_{2}^{'}+\mathbf{R}\mathbf{A}_{1}^{'}+\mathbf{C}_{1}^{'}.
\end{equation}

\begin{lemma} \label{lem:4.1}
If $\widetilde{\alpha}<1$, the equation $(4.1)$ has the minimal non-negative solution
\begin{equation}\label{4.2}
    \mathbf{R}=\left[\begin{array}{cc}
   r & \frac{r\theta}{\bar{\theta}\overline{p(1)}\mu_{b}(1-r)} \\
   0 & \widetilde{\alpha}
   \end{array}\right]
\end{equation}
where
\begin{equation*}
r=\frac{1}{2\mu_{\nu}\overline{p(0)}}\left[\beta+p(0)\bar{\mu}_{\nu}+\overline{p(0)}\mu_{\nu}-\sqrt{(\beta+p(0)\bar{\mu}_{\nu}+\overline{p(0)}\mu_{\nu})^{2}-4p(0)\mu_{\nu}\overline{p(0)}\bar{\mu}_{\nu}}\right],
~\beta=\theta{\bar{\theta}}^{-1},
\end{equation*}
and $0<r<1$.
\end{lemma}

\proof As $\mathbf{A}_{1}^{'}$, $\mathbf{C}_{1}^{'}$, $\mathbf{B}_{2}^{'}$ are all upper-triangular matrices, $\mathbf{R}$ has the same form. Then we can suppose that
\[
    \mathbf{R}=\left [
 \begin{array}{cc}
   r_{11} & r_{12} \\
   0 & r_{22}
   \end{array}\right ].
   \]
Substituting $\mathbf{A}_{1}^{'}$, $\mathbf{C}_{1}^{'}$, $\mathbf{B}_{2}^{'}$ into (\ref{4.1}), we get the equations as follows:
\begin{equation}\label{4.3}
     \left\{
     \begin{aligned}
                    r_{11}&=\overline{p(0)}\mu_{\nu}\bar{\theta}r_{11}^{2}+\bar{\theta}(1-p(0)\bar{\mu}_{\nu}-\overline{p(0)}\mu_{\nu})r_{11}+p(0)\bar{\theta}\bar{\mu}_{\nu},\\
                    r_{12}&=\overline{p(0)}\theta\mu_{\nu}r_{11}^{2}+\overline{p(1)}\mu_{b}r_{11}r_{12}+\overline{p(1)}\mu_{b}r_{12}r_{22}+\theta(1-p(0)\bar{\mu}_{\nu}-\overline{p(0)}\mu_{\nu})r_{11}\\
                          &~~~+(1-p(1)\bar{\mu}_{b}-\overline{p(1)}\mu_{b})r_{12}+p(0)\theta\bar{\mu}_{\nu},\\
                    r_{22}&=\overline{p(1)}\mu_{b}r_{22}^{2}+(1-p(1)\bar{\mu}_{b}-\overline{p(1)}\mu_{b})r_{22}+p(1)\bar{\mu}_{b}.
                  \end{aligned}
\right.
\end{equation}
The third equation of (\ref{4.3}) has the minimal non-negative solution $r_{22}=\widetilde{\alpha}$ (the other root is $r_{22}=1$). The first equation of (\ref{4.3}) can be written as:
\begin{equation*}
\overline{p(0)}\mu_{\nu}r_{11}^{2}-(\beta+p(0)\bar{\mu}_{\nu}+\overline{p(0)}\mu_{\nu})r_{11}+p(0)\bar{\mu}_{\nu}=0.
\end{equation*}
The root of the above equation is $r_{11}=r$, and $0<r<1$ (the other root is greater than 1). Substituting $r_{11}=r$ and $r_{22}=\widetilde{\alpha}$ into the second equation, we can get
\begin{equation*}
r_{12}=\frac{r\theta}{{\bar{\theta}}\overline{p(1)}\mu_{b}(1-r)}.
\end{equation*} \BOX\\
Besides, $r$ satisfies $\frac{r\theta}{\overline{\theta}(1-r)}=p(0)\bar{\mu}_{\nu}-r\overline{p(0)}\mu_{\nu}$.

\begin{theorem} \label{the:4.1}
Consider a partially observable $Geo/Geo/1$ queue with multiple working vocations, in which all arriving customers follow the mixed policy $(q(0),q(1))$. If $\widetilde{\alpha}<1$, the stationary distribution $\{\pi_{kj}^{+'}|(k,j)\in\Omega_{po}\}$ of $(L^{+},J)$ is
\begin{equation}\label{4.4}
     \left\{
     \begin{aligned}
                    \pi_{00}^{+'}&=K[\theta+\bar{\theta}\overline{p(0)}\mu_{\nu}(1-r)],\\
                    \pi_{k0}^{+'}&=Kp(0)\bar{\theta}(1-r)r^{k-1}, ~k\geq1,\\
                    \pi_{k1}^{+'}&=K\frac{p(0)\theta}{\overline{p(1)}\mu_{b}}\sum_{j=0}^{k-1}r^{j}\widetilde{\alpha}^{k-1-j}, ~k\geq1,
                  \end{aligned}
\right.
\end{equation}
where
\begin{equation*}
K=\frac{\overline{p(1)}\mu_{b}(1-r)(1-\widetilde{\alpha})}{\overline{p(1)}\mu_{b}(1-r)(1-\widetilde{\alpha})[\theta+\bar{\theta}\overline{p(0)}\mu_{\nu}(1-r)+p(0)\bar{\theta}]+p(0)\theta}.
\end{equation*}
\end{theorem}

\proof According to the matrix-geometric solution method in Tian et al. \cite{20}, we have
\begin{equation}\label{4.5}
\pi_{k}^{'}=(\pi_{k0}^{+'},\pi_{k1}^{+'})=(\pi_{10}^{+'},\pi_{11}^{+'})\mathbf{R}^{k-1}, \ \ k\geq1,
\end{equation}
and $(\pi_{00}^{+'},\pi_{10}^{+'},\pi_{11}^{+'})$ satisfies the equations
\begin{equation*}
(\pi_{00}^{+'},\pi_{10}^{+'},\pi_{11}^{+'})\mathbf{B}[\mathbf{R}]=(\pi_{00}^{+'},\pi_{10}^{+'},\pi_{11}^{+'}).
\end{equation*}
where
\begin{equation*}
\begin{aligned}
\mathbf{B}[\mathbf{R}]&=\left [
\begin{array}{cc}
\mathbf{B}_{0}^{'} & \mathbf{C}_{0}^{'} \\
\mathbf{B}_{1}^{'} & \mathbf{R}\mathbf{B}_{2}^{'}+\mathbf{A}_{1}^{'}
\end{array}\right ]\\
&=\left [
\begin{array}{ccc}
\overline{p(0)} & p(0)\bar{\theta} & p(0)\theta\\
\overline{p(0)}\mu_{\nu} & \bar{\theta}(1-\overline{p(0)}\mu_{\nu})-\frac{r\theta}{1-r} & \theta(1-\overline{p(0)}\mu_{\nu})+\frac{r\theta}{1-r}\\
\overline{p(1)}\mu_{b} & 0 & 1-\overline{p(1)}\mu_{b}
\end{array}\right ].
\end{aligned}
\end{equation*}
Solving the above equations, we get
\begin{equation*}
     \left\{
     \begin{aligned}
                    \pi_{00}^{+'}&=\overline{p(0)}\pi_{00}^{+'}+\overline{p(0)}\mu_{\nu}\pi_{10}^{+'}+\overline{p(1)}\mu_{b}\pi_{11}^{+'},\\
                    \pi_{10}^{+'}&=p(0)\bar{\theta}\pi_{00}^{+'}+\left[\bar{\theta}(1-\overline{p(0)}\mu_{\nu})-\frac{r\theta}{1-r}\right]\pi_{10}^{+'},\\
                    \pi_{11}^{+'}&=p(0)\theta\pi_{00}^{+'}+\left[\theta(1-\overline{p(0)}\mu_{\nu})+\frac{r\theta}{1-r}\right]\pi_{10}^{+'}+(1-\overline{p(1)}\mu_{b})\pi_{11}^{+'}.
                  \end{aligned}
\right.
\end{equation*}
Setting $\pi_{0}^{'}=\pi_{00}^{+'}$ is a constant which can be determined by the normalization condition $\pi_{0}^{'}\bm{e}+\pi_{1}^{'}(\mathbf{I}-\mathbf{R})^{-1}\bm{e}=1$.\\
Hence,
\begin{equation*}
\pi_{10}^{+'}=\frac{p(0)\bar{\theta}(1-r)}{\theta+\bar{\theta}\overline{p(0)}\mu_{\nu}(1-r)}\pi_{00}^{+'},
\end{equation*}
\begin{equation*}
\pi_{11}^{+'}=\frac{p(0)\theta}{\overline{p(1)}\mu_{b}[\theta+\bar{\theta}\overline{p(0)}\mu_{\nu}(1-r)]}\pi_{00}^{+'}.
\end{equation*}
From (\ref{4.2}), we have
\begin{equation*}
 \mathbf{R}^{k}=\left [
 \begin{array}{cc}
   r^{k} & \frac{r\theta}{\bar{\theta}\overline{p(1)}\mu_{b}(1-r)}\sum_{j=0}^{k-1}r^{j}\widetilde{\alpha}^{k-1-j} \\
   0 & \widetilde{\alpha}^{k}
   \end{array}\right ],
~~k\geq1.
\end{equation*}
Substituting $\pi_{1}^{'}=(\pi_{10}^{+'},\pi_{11}^{+'})$ and $\mathbf{R}^{k-1}$ into (\ref{4.5}), we obtain
\begin{equation*}
\begin{aligned}
\pi_{k}^{'}&=(\pi_{k0}^{+'},\pi_{k1}^{+'})\\
&=\frac{\pi_{00}^{+'}}{\theta+\bar{\theta}\overline{p(0)}\mu_{\nu}(1-r)}\left(p(0)\bar{\theta}(1-r)r^{k-1},\frac{p(0)\theta}{\overline{p(1)}\mu_{b}}\sum_{j=0}^{k-1}r^{j}\widetilde{\alpha}^{k-1-j}\right),~~k\geq1.
\end{aligned}
\end{equation*}
In the end, using the normalization condition
\begin{equation*}
\pi_{0}^{'}\bm{e}+\pi_{1}^{'}(\mathbf{I}-\mathbf{R})^{-1}\bm{e}=1,
\end{equation*}
we get
\begin{equation*}
\pi_{00}^{+'}=\frac{[\theta+\bar{\theta}\overline{p(0)}\mu_{\nu}(1-r)]\overline{p(1)}\mu_{b}(1-r)(1-\widetilde{\alpha})}{\overline{p(1)}\mu_{b}(1-r)(1-\widetilde{\alpha})[\theta+\bar{\theta}\overline{p(0)}\mu_{\nu}(1-r)+p(0)\bar{\theta}]+p(0)\theta}.
\end{equation*}
So we complete the proof.   \BOX

According to (\ref{4.4}), we can calculate the probabilities that an arrival occurs in a regular busy period and in a working vacation period respectively as follows:
\begin{equation}\label{4.6}
\begin{aligned}
P(J=0)=\sum_{k=0}^{\infty}\pi_{k0}^{+'}
=K\left(\frac{\theta}{1-r}+\bar{\theta}\mu_{\nu}\right),
\end{aligned}
\end{equation}
\begin{equation}\label{4.7}
\begin{aligned}
P(J=1)=\sum_{k=1}^{\infty}\pi_{k1}^{+'}
=\frac{Kp(0)\theta}{\overline{p(1)}\mu_{b}(1-\widetilde{\alpha})(1-r)}.
\end{aligned}
\end{equation}
Then we can compute the conditional mean queue length in a working vacation period and in a regular busy period respectively as follows:
\begin{equation*}
E[L_{0}]=\frac{\sum_{k=1}^{\infty}k\pi_{k0}^{+'}}{\sum_{k=0}^{\infty}\pi_{k0}^{+'}}
=\frac{p(0)\bar{\theta}}{\theta+\bar{\theta}\mu_{\nu}(1-r)},
\end{equation*}
\begin{equation*}
E[L_{1}]=\frac{\sum_{k=1}^{\infty}k\pi_{k1}^{+'}}{\sum_{k=1}^{\infty}\pi_{k1}^{+'}}=\frac{1-r\widetilde{\alpha}}{(1-\widetilde{\alpha})(1-r)}.
\end{equation*}

As for the conditional mean sojourn time in regular busy period, we can similarly consider two situations at state $(k,1)$ $(k\geq1)$ with those in observable queues. So we get the PGF of the customers' sojourn time in regular busy period, denoted by $W_{1}^{*}(z)$.
\begin{equation*}
\begin{aligned}
W_{1}^{*}(z)&=\sum_{k=1}^{\infty}\pi_{k1}^{+'}\left[\mu_{b}\left(\frac{\mu_{b}z}{1-\bar{\mu}_{b}z}\right)^{k}+\bar{\mu}_{b}\left(\frac{\mu_{b}z}{1-\bar{\mu}_{b}z}\right)^{k+1}\right]\\
&=\frac{Kp(0)\theta\mu_{b}z}{\overline{p(1)}(1-\bar{\mu}_{b}z-\mu_{b}\widetilde{\alpha} z)(1-\bar{\mu}_{b}z-\mu_{b}rz)}.
\end{aligned}
\end{equation*}
Hence, the PGF of the conditional sojourn time in regular busy period is
\begin{equation*}
{W_{1}}(z)=\frac{W_{1}^{*}(z)}{P(J=1)}=\frac{{\mu_{b}}^{2}z(1-\widetilde{\alpha})(1-r)}{(1-\bar{\mu}_{b}z-\mu_{b}\widetilde{\alpha}z)(1-\bar{\mu}_{b}z-\mu_{b}rz)}.
\end{equation*}
Then we obtain the conditional mean sojourn time of the customers in regular busy period as \begin{equation*}
E[W_{1}]=W_{1}^{'}(z)|_{z=1}=\frac{1}{\mu_{b}(1-r)}+\frac{\mu_{b}\widetilde{\alpha}-\mu_{b}+1}{\mu_{b}(1-\widetilde{\alpha})}.
\end{equation*}

As for the conditional mean sojourn time in working vacation period, we can similarly consider two situations at state $(k,0)$ $(k\geq0)$ with those in observable queues. So we get the PGF of the customers' sojourn time in working vacation period, denoted by $W_{0}^{*}(z)$.
\begin{equation*}
\begin{aligned}
W_{0}^{*}(z)&=\sum_{k=1}^{\infty}\pi_{k0}^{+'}\mu_{\nu}\left[P(S_{\nu}^{(k)}\leq V)S_{\nu}^{(k)}(z|S_{\nu}^{(k)}\leq V)\right.\\
&~~~\left.+\sum_{j=0}^{k-1}P(S_{\nu}^{(j)}\leq V <S_{\nu}^{(j+1)})V(z|S_{\nu}^{(j)}\leq V <S_{\nu}^{(j+1)})\left(\frac{\mu_{b}}{1-\bar{\mu}_{b}z}\right)^{k-j}\right]\\
&~~~+\sum_{k=1}^{\infty}\pi_{k0}^{+'}\bar{\mu}_{\nu}\left[P(S_{\nu}^{(k+1)}\leq V)S_{\nu}^{(k+1)}(z|S_{\nu}^{(k+1)}\leq V)\right.\\
&~~~\left.+\sum_{j=0}^{k}P(S_{\nu}^{(j)}\leq V <S_{\nu}^{(j+1)})V(z|S_{\nu}^{(j)}\leq V <S_{\nu}^{(j+1)})\left(\frac{\mu_{b}}{1-\bar{\mu}_{b}z}\right)^{k+1-j}\right]\\
&~~~+\pi_{00}^{+'}\left[P(S_{\nu}^{(1)}\leq V)S_{\nu}^{(1)}(z|S_{\nu}^{(1)}\leq V)+P(S_{\nu}^{(0)}\leq V <S_{\nu}^{(1)})V(z|S_{\nu}^{(0)}\leq V <S_{\nu}^{(1)})\right.\\
&~~~\left.\cdot\frac{\mu_{b}}{1-\bar{\mu}_{b}z}\right]\\
&=\sum_{k=1}^{\infty}\pi_{k0}^{+'}\mu_{\nu}\left[\left(\frac{\mu_{\nu}\bar{\theta}z}{1-\bar{\mu}_{\nu}\bar{\theta}z}\right)^{k}
+\sum_{j=0}^{k-1}\frac{\theta}{1-\bar{\mu}_{\nu}\bar{\theta}z}\left(\frac{\mu_{\nu}\bar{\theta}z}{1-\bar{\mu}_{\nu}\bar{\theta}z}\right)^{j}\left(\frac{\mu_{b}z}{1-\bar{\mu}_{b}z}\right)^{k-j}\right]\\
&~~~+\sum_{k=1}^{\infty}\pi_{k0}^{+'}\bar{\mu}_{\nu}\left[\left(\frac{\mu_{\nu}\bar{\theta}z}{1-\bar{\mu}_{\nu}\bar{\theta}z}\right)^{k+1}
+\sum_{j=0}^{k}\frac{\theta}{1-\bar{\mu}_{\nu}\bar{\theta}z}\left(\frac{\mu_{\nu}\bar{\theta}z}{1-\bar{\mu}_{\nu}\bar{\theta}z}\right)^{j}\left(\frac{\mu_{b}z}{1-\bar{\mu}_{b}z}\right)^{k+1-j}\right]\\
&~~~+\pi_{00}^{+'}\left(\frac{\mu_{\nu}\bar{\theta}z}{1-\bar{\mu}_{\nu}\bar{\theta}z}+\frac{\theta}{1-\bar{\mu}_{\nu}\bar{\theta}z}\frac{\mu_{b}z}{1-\bar{\mu}_{b}z}\right)\\
&=\frac{Kp(0)(1-r)\mu_{\nu}^{2}\bar{\theta}^{2}z}{(1-\bar{\mu}_{\nu}\bar{\theta}z)(1-\bar{\mu}_{\nu}\bar{\theta}z-\mu_{\nu}\bar{\theta}rz)}
+\frac{Kp(0)(1-r)\theta\bar{\theta}\mu_{b}z}{(1-\bar{\mu}_{b}z-\mu_{b}rz)(1-\bar{\mu}_{\nu}\bar{\theta}z-\mu_{\nu}\bar{\theta}rz)}\left(\mu_{\nu}+\frac{\bar{\mu}_{\nu}}{r}\right)\\
&~~~+K[\theta+\bar{\theta}\overline{p(0)}\mu_{\nu}(1-r)]\left(\frac{\mu_{\nu}\bar{\theta}z}{1-\bar{\mu}_{\nu}\bar{\theta}z}+\frac{\theta}{1-\bar{\mu}_{\nu}\bar{\theta}z}\frac{\mu_{b}z}{1-\bar{\mu}_{b}z}\right).
\end{aligned}
\end{equation*}
Hence, the PGF of the conditional sojourn time in working vacation period is
\begin{equation*}
\begin{aligned}
W_{0}(z)&=\frac{W_{0}^{*}(z)}{P(J=0)}\\
&=\frac{1}{\frac{\theta}{(1-r)}+\bar{\theta}\mu_{\nu}}\left\{\frac{p(0)(1-r)\mu_{\nu}^{2}\bar{\theta}^{2}z}{(1-\bar{\mu}_{\nu}\bar{\theta}z)(1-\bar{\mu}_{\nu}\bar{\theta}z-\mu_{\nu}\bar{\theta}rz)}
+\frac{p(0)(1-r)\theta\bar{\theta}\mu_{b}z}{(1-\bar{\mu}_{b}z-\mu_{b}rz)(1-\bar{\mu}_{\nu}\bar{\theta}z-\mu_{\nu}\bar{\theta}rz)}\right.\\
&~~~\left.\cdot\left(\mu_{\nu}+\frac{\bar{\mu}_{\nu}}{r}\right)+[\theta+\bar{\theta}\overline{p(0)}\mu_{\nu}(1-r)]
\left(\frac{\mu_{\nu}\bar{\theta}z}{1-\bar{\mu}_{\nu}\bar{\theta}z}
+\frac{\theta}{1-\bar{\mu}_{\nu}\bar{\theta}z}\frac{\mu_{b}z}{1-\bar{\mu}_{b}z}\right)\right\}.
\end{aligned}
\end{equation*}
Then we obtain the conditional mean sojourn time of the customers in working vacation period as
\begin{equation*}
\begin{aligned}
E[W_{0}]&=W_{0}^{'}(z)|_{z=1}\\
&=\frac{p(0)(1-r)^{2}\mu_{\nu}^{2}\bar{\theta}^{2}[\mu_{\nu}(2-r)-\mu_{\nu}^{2}(1-r)+\theta(2-\theta)(1-\mu_{\nu})(1+r\mu_{\nu}-\mu_{\nu})]}{(\theta+\mu_{\nu}-\theta\mu_{\nu})^{2}(\theta+\mu_{\nu}-\theta\mu_{\nu}-r\mu_{\nu}+r\theta\mu_{\nu})^{3}}\\
&~~~+\frac{p(0)\bar{\theta}\theta[\mu_{\nu}(1-r)+\mu_{b}(1-r)(1+r\mu_{\nu}-\mu_{\nu})+\theta(1+r\mu_{\nu}-\mu_{\nu})(1+r\mu_{b}-\mu_{b})]}{\mu_{b}(\theta+\mu_{\nu}-\theta\mu_{\nu}-r\mu_{\nu}+r\theta\mu_{\nu})^{3}}\\
&~~~+\frac{[\theta+\overline{p(0)}\bar{\theta}\mu_{\nu}(1-r)](\theta+\mu_{b}-\theta\mu_{b})(1-r)}{\mu_{b}(\theta+\mu_{\nu}-\theta\mu_{\nu})(\theta+\mu_{\nu}-\mu_{\nu}r-\mu_{\nu}\theta+r\mu_{\nu}\theta)}.
\end{aligned}
\end{equation*}

Hence, the expected net benefit of a customer who joins the queueing system equals
\begin{equation*}
\begin{aligned}
U(0;q(0))&=R-C\left\{\frac{p(0)(1-r)^{2}\mu_{\nu}^{2}\bar{\theta}^{2}[\mu_{\nu}(2-r)-\mu_{\nu}^{2}(1-r)+\theta(2-\theta)(1-\mu_{\nu})(1+r\mu_{\nu}-\mu_{\nu})]}{(\theta+\mu_{\nu}-\theta\mu_{\nu})^{2}(\theta+\mu_{\nu}-\theta\mu_{\nu}-r\mu_{\nu}+r\theta\mu_{\nu})^{3}}\right.\\
&~~~+\frac{p(0)\bar{\theta}\theta[\mu_{\nu}(1-r)+\mu_{b}(1-r)(1+r\mu_{\nu}-\mu_{\nu})+\theta(1+r\mu_{\nu}-\mu_{\nu})(1+r\mu_{b}-\mu_{b})]}{\mu_{b}(\theta+\mu_{\nu}-\theta\mu_{\nu}-r\mu_{\nu}+r\theta\mu_{\nu})^{3}}\\
&~~~\left.+\frac{[\theta+\overline{p(0)}\bar{\theta}\mu_{\nu}(1-r)](\theta+\mu_{b}-\theta\mu_{b})(1-r)}{\mu_{b}(\theta+\mu_{\nu}-\theta\mu_{\nu})(\theta+\mu_{\nu}-\mu_{\nu}r-\mu_{\nu}\theta+r\mu_{\nu}\theta)}\right\},
\end{aligned}
\end{equation*}
\begin{equation*}
U(1;q(0),q(1))=R-C\left[\frac{1}{\mu_{b}(1-r)}+\frac{\mu_{b}\widetilde{\alpha}-\mu_{b}+1}{\mu_{b}(1-\widetilde{\alpha})}\right].
\end{equation*}
Solving $U(0;q(0))=0$ and $U(1;q(0),q(1))=0$, we
obtain the positive and feasible roots $(q_{e}^{*}(0),q_{e}^{*}(1))$.
Because the uniqueness of the roots in the partially observable case,
we present a set of numerical experiments to show the effect of the information level.
We get that the customers' equilibrium mixed strategy $(q_{e}(0),q_{e}(1))$ is
 unique when $(q_{e}(0),q_{e}(1))=(min\{q_{e}^{*}(0),1\},min\{q_{e}^{*}(1),1\})$ in the partially observable case.
\begin{figure}
\centering
\includegraphics[width=10cm,height=8cm]{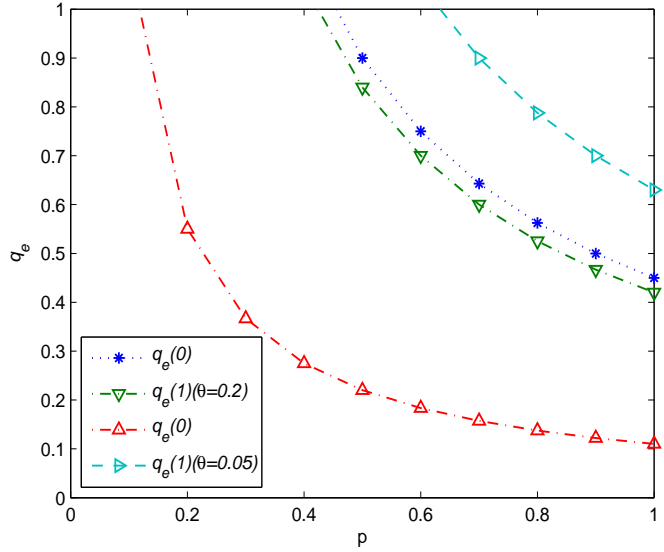}
\caption{Equilibrium thresholds for the partially observable systems with different working vacation parameter $\theta$ when $R=8$, $C=3$, $\mu_{b}=0.8$, $\mu_{\nu}=0.4$}\label{Fig.6}
\end{figure}

From the Figure~\ref{Fig.6} we can obtain that despite the server¡¯s
service rates have the relation $\mu_{b}>\mu_{\nu}$,
Figure~\ref{Fig.6} shows that $q_{e}(0)>q_{e}(1)$. The reason for
this phenomenon is that the lower service rate $\mu_{\nu}$ in
working vacation state can meet the customers' service demand when
the waiting cost is low and the expected working vacation period is
relatively short.

Next, we consider the equilibrium and maximal social benefit. Using theorem~\ref{the:4.1}, we can get the mean queue length
\begin{equation*}
\begin{aligned}
E[L]=\sum_{k=0}^{\infty}k\pi_{k0}^{+'}+\sum_{k=1}^{\infty}k\pi_{k1}^{+'}
=\frac{Kp(0)}{1-r}\left[\bar{\theta}+\frac{\theta(1-r\widetilde{\alpha})}{\overline{p(1)}\mu_{b}(1-r)(1-\widetilde{\alpha})^{2}}\right].
\end{aligned}
\end{equation*}
So the social benefit per time unit for the mixed policy $(q(0),q(1))$ can be calculated as
\begin{equation*}
\begin{aligned}
U_{s}(q(0),q(1))&=dR-CE[L]\\
&=p(P(J=0)q(0)+P(J=1)q(1))R-\frac{CKp(0)}{1-r}\left[\bar{\theta}+\frac{\theta(1-r\widetilde{\alpha})}{\overline{p(1)}\mu_{b}(1-r)(1-\widetilde{\alpha})^{2}}\right].
\end{aligned}
\end{equation*}
where $d=p(P(J=0)q(0)+P(J=1)q(1))$ and $P(J=0)$, $P(J=1)$ are given by (\ref{4.6}) and (\ref{4.7}).

\begin{figure}
\centering
\includegraphics[width=10cm,height=8cm]{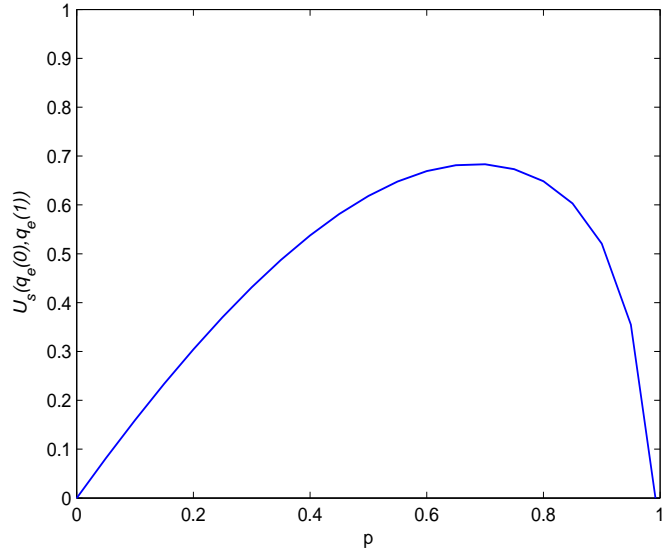}
\caption{Equilibrium social benefit for the partially observable systems with $R=10$, $C=3$, $\theta=0.05$, $\mu_{\nu}=0.5$, $\mu_{b}=0.9$}\label{Fig.7}
\end{figure}

When all customers follow the above equilibrium mixed strategy
$(q_{e}(0),q_{e}(1))$, the social benefit per time unit in
equilibrium can be represented as $U_{s}(q_{e}(0),q_{e}(1))$.
Figure~\ref{Fig.7} is concerned with the social benefit under the
equilibrium mixed strategy. We can observe that
$U_{s}(q_{e}(0),q_{e}(1))$ first increases, then decreases with
respect to $p$.

\begin{figure}
\centering
\includegraphics[width=10cm,height=8cm]{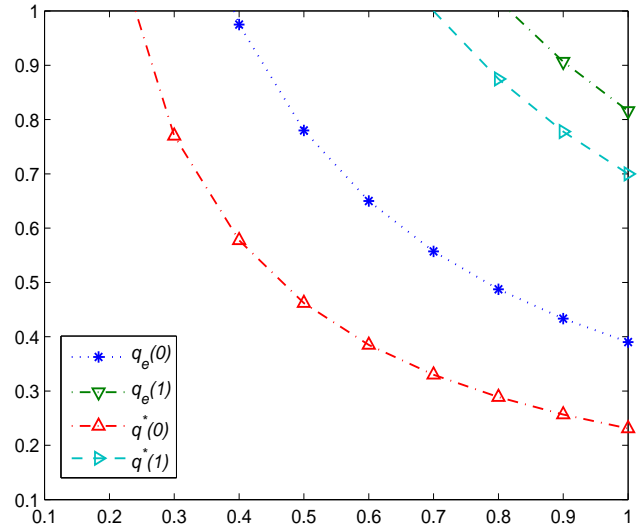}
\caption{Comparisons of equilibrium and socially optimal mixed strategies for the partially observable systems with $R=10$, $C=3$, $\theta=0.05$, $\mu_{\nu}=0.5$, $\mu_{b}=0.9$}\label{Fig.8}
\end{figure}

Next, from the view of social optimization, denote the socially
optimal mixed strategy as $(q^{*}(0),q^{*}(1))$, which can be
obtained by differentiating $U_{s}(q(0),q(1))$ with respect to
$q(0)$ and $q(1)$ and solving the equations. Figure~\ref{Fig.8}
compares the equilibrium mixed strategy $(q_{e}(0),q_{e}(1))$
 and the socially optimal mixed strategy $(q^{*}(0),q^{*}(1))$ for the partially observable systems.
  We get that $q_{e}(0)>q^{*}(0)$ and $q_{e}(1)>q^{*}(1)$,
  which shows that the individual optimization deviates from the social optimization
   for the partially observable systems.

\section{The unobservable queue}
\label{5}

\begin{figure}
\centering
\includegraphics[width=15cm,height=5.5cm]{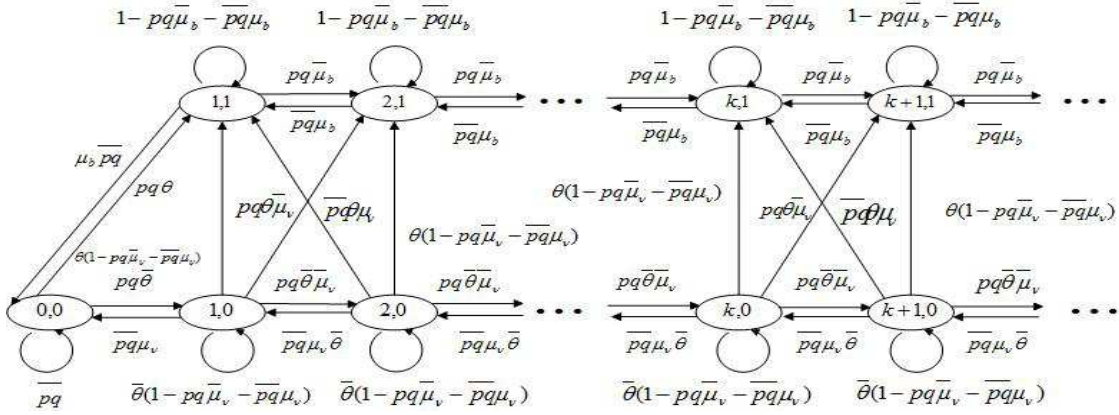}
\caption{Transition rate diagram for the unobservable queues}\label{Fig.9}
\end{figure}

We finally investigate the unobservable case in which the arriving customers observe
neither the state $J_{n}$ of the server nor the number of customers $L_{n}^{+}$ in the system.
 Suppose that the customers' decision that whether to join or to balk upon their arrival can be indicated
 by a joining probability $q~(0\leq q\leq1)$, and their equilibrium mixed strategy is represented by $q_{e}$.
  The state space is $\Omega_{uo}=\{(0,0)\}\bigcup\{(k,j):k\geqslant1,~j=0,1\}$. The transition rate
diagram is described in Figure~\ref{Fig.9}.

According to the decomposition result given in Tian et al. \cite{20}, we can obtain the mean sojourn time of a joining customer as follows:
\begin{equation*}
E[W]=\frac{1}{\mu_{b}(1-\alpha^{'})}+K^{*}\delta_{2}(pq)^{-1}(1-\overline{pq}\sigma)\frac{\sigma}{1-\sigma},
\end{equation*}
where
\begin{eqnarray*}
K^{*}=[\delta_{1}(pq)^{-1}+\delta_{2}(pq)^{-1}(1-\overline{pq}\sigma)]^{-1},
\end{eqnarray*}
\[
\delta_{1}=(pq)^{2}{\bar{\mu}}_{b}\bar{\theta}\mu_{\nu}\frac{(1-r^{'})^{2}}{r^{'}}, \ \
\delta_{2}=pq\bar{\theta}(pq+r^{'}\overline{pq})\frac{1-r^{'}}{r^{'}}(\mu_{b}-\mu_{\nu}),
\]
\[
\alpha^{'}=\frac{pq{\bar{\mu}}_{b}}{\overline{pq}\mu_{b}}, \ \ \sigma=\frac{r^{'}}{pq+r^{'}\overline{pq}}, \ \ 1-\sigma=\frac{pq(1-r^{'})}{pq+r^{'}\overline{pq}},\\
\]
and $r^{'}$ satisfies the equation $pq{\bar{\mu}}_{\nu}-r^{'}\mu_{\nu}\overline{pq}=\frac{r^{'}\theta}{\overline{\theta}(1-r^{'})}$.

\begin{figure}
\centering
\includegraphics[width=10cm,height=8cm]{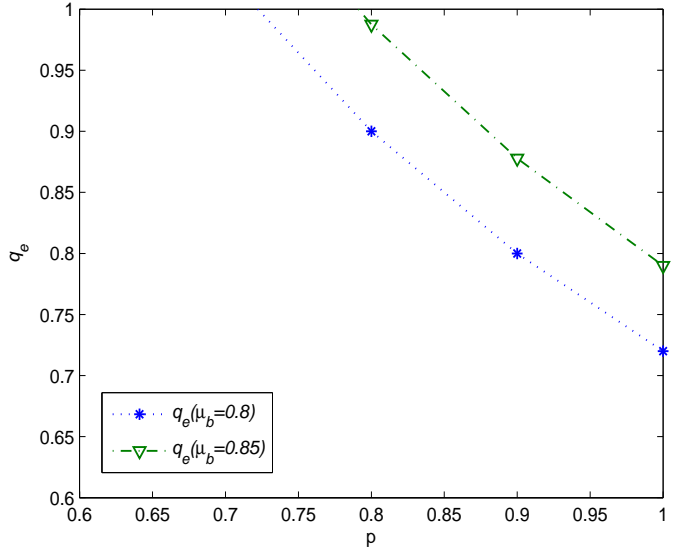}
\caption{Equilibrium thresholds for the unobservable systems with different service parameter $\mu_{b}$ when $R=4.5$, $C=1$, $\theta=0.3$, $\mu_{\nu}=0.5$}\label{Fig.10}
\end{figure}

Hence, the expected net benefit of the customer who joins the queueing system equals
\begin{equation*}
\begin{aligned}
U(q)&=R-CE[W]\\
&=R-C\left[\frac{1}{\mu_{b}(1-\alpha^{'})}+K^{*}\delta_{2}(pq)^{-1}(1-\overline{pq}\sigma)\frac{\sigma}{1-\sigma}\right].
\end{aligned}
\end{equation*}
Solving $U(q)=0$, we obtain the positive and feasible root $q_{e}^{*}$.
 Because the uniqueness of the root in the unobservable case,
 we present a set of numerical experiments to show the effect of the information level.
  We get that the customers' equilibrium mixed strategy $q_{e}$ is unique when $q_{e}=min\{q_{e}^{*},1\}$
  in the unobservable case. From the Figure~\ref{Fig.10} we can observe that
$q_{e}$ decreases with respect to $p$, and $q_{e}$ decreases faster
when $\mu_{b}$ becomes smaller.

Next, we can get the social benefit per time unit for the mixed policy $q$ as
\begin{equation*}
\begin{aligned}
U_{s}(q)&=pq(R-CE[W])\\
&=pq\left[R-C\left(\frac{1}{\mu_{b}(1-\alpha^{'})}+K^{*}\delta_{2}(pq)^{-1}(1-\overline{pq}\sigma)\frac{\sigma}{1-\sigma}\right)\right].
\end{aligned}
\end{equation*}

\begin{figure}
\centering
\includegraphics[width=10cm,height=8cm]{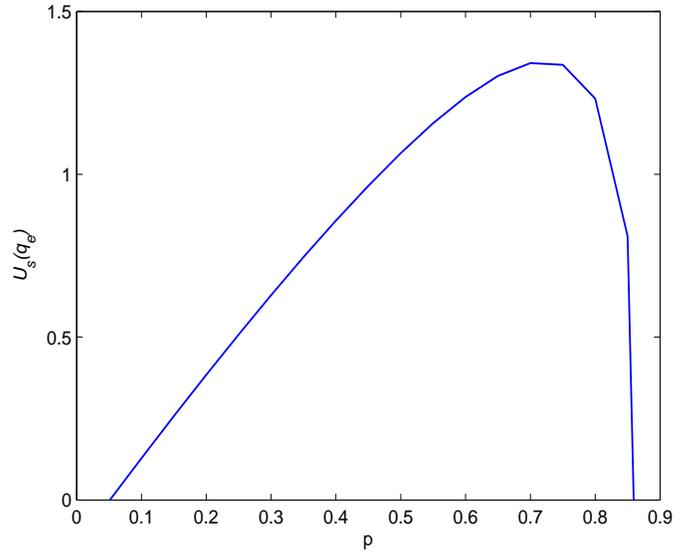}
\caption{Equilibrium social benefit for the unobservable systems with $R=4.5$, $C=1$, $\theta=0.3$, $\mu_{\nu}=0.5$, $\mu_{b}=0.9$}\label{Fig.11}
\end{figure}

When all customers follow the above equilibrium mixed strategy
$q_{e}$, the social benefit per time unit in equilibrium can be
represented as $U_{s}(q_{e})$. Figure~\ref{Fig.11} is concerned with
the social benefit under the equilibrium mixed strategy. We can
observe that $U_{s}(q_{e})$ first increases, then decreases with
respect to $p$.

\begin{figure}
\centering
\includegraphics[width=10cm,height=8cm]{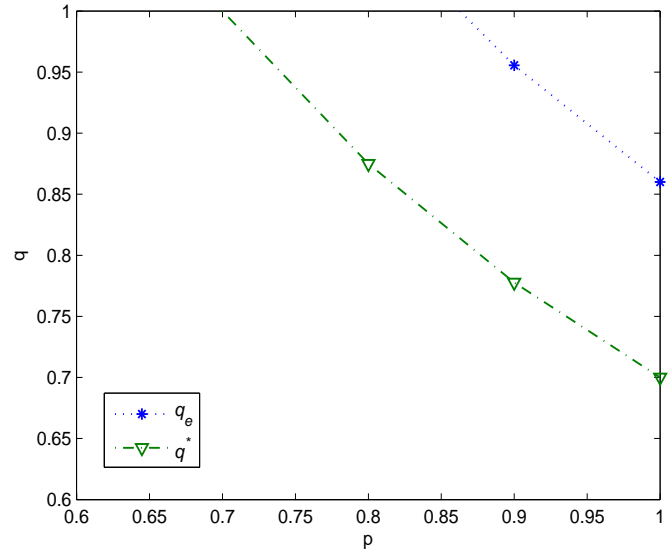}
\caption{Comparisons of equilibrium and socially optimal mixed strategies for the unobservable systems with $R=4.5$, $C=1$, $\theta=0.3$, $\mu_{\nu}=0.5$, $\mu_{b}=0.9$}\label{Fig.12}
\end{figure}

Next, from the view of social optimization, denote the socially
optimal mixed strategy as $q^{*}$ , which can be obtained by
differentiating $U_{s}(q)$ with respect to $q$ and solving the
equation. Figure~\ref{Fig.12} compares the equilibrium mixed
strategy $q_{e}$ and the socially optimal mixed strategy $q^{*}$ for
the unobservable systems. We get that $q_{e}>q^{*}$, which shows
that the individual optimization deviates from the social
optimization for the unobservable systems.

\section{Conclusions}
\label{6}

In this paper we have researched the equilibrium customer behavior in the discrete-time $Geo/Geo/1$ queueing system with multiple working vacations in which arriving customers decide whether to join or to balk the system. To the best of the authors' knowledge, it is the first time studying the discrete-time queueing systems with multiple working vacations from an economic viewpoint. Three different cases with respect to the level of information provided to arriving customers have been investigated extensively and the equilibrium thresholds and social benefit for each case were derived. Furthermore, we have presented some numerical experiments to research the effect of the information level on the equilibrium behavior and to compare the customers' equilibrium and socially optimal strategies.

The focal point of this paper is on the equilibrium balking strategies' analysis. There are various aspects for future research. One can think of the equilibrium behavior under single working vacation policy. Furthermore, we can also explore equilibrium behavior in $Geo/G/1$ and $GI/Geo/1$ queueing systems with various vacation policies.

\section*{Acknowledgements}
This research is partially supported by the National Natural Science
Foundation of China (11201489, 11271373, 11371374).

\end{document}